# Generalized weight properties of resultants and discriminants, and applications to projective enumerative geometry

Laurent Busé, Thomas Dedieu

**Abstract.** In a book dating back to 1862, Salmon stated a formula giving the first terms of the Taylor expansion of the discriminant of a plane algebraic curve, and from it derived various enumerative quantities for surfaces in $\mathbf{P}^3$. In this text, we provide complete proofs of this formula and its enumerative applications, and extend Salmon's considerations to hypersurfaces in a projective space of arbitrary dimension. To this end, we extend reduced elimination theory by introducing the concept of reduced discriminant, and provide a thorough study of its weight properties; the latter are deeply linked to projective enumerative geometric properties. Then, following Salmon's approach, we compute the number of members of a pencil of hyperplanes that are bitangent to a fixed projective hypersurface. Some other results in the same spirit are also discussed.

In his book [16] originally published in 1862, Salmon casually gives the first terms of the Taylor expansion of the discriminant of a plane algebraic curve of degree $d \geqslant 2$. In a suitable system of homogeneous coordinates, any plane curve $V(f) \subset \mathbf{P}^2$ has an equation of the form

$$f(x,y,z) = Tz^{d-1}y + \frac{1}{2}z^{d-2}(Ax^2 + 2Bxy + Cy^2) + \sum_{k \geqslant 3} z^{d-k} f_k(x,y) = 0$$

where each polynomial $f_k(x,y)$ is homogeneous of degree $k$ in $x, y$. In this notation, Salmon states as a well known fact and without proof that the discriminant of $f(x, y, z)$ has the form

(🦉) $$\operatorname{Disc}(f) = T^2 A(B^2 - AC)^2 \varphi + T^3 \psi$$

where $\varphi$ is "the discriminant when $T$ vanishes" [16, § 605]. This note arose as an attempt to prove this formula and to shed light on the geometric content of the vanishing of $\varphi$. It turns out that this polynomial is deeply linked to the concept of *reduced resultant* introduced by Zariski much later in 1936 [18]; this leads us to introduce the *reduced discriminant* of a hypersurface, of which $\varphi$ is an instance. Of course when $T$ is zero, $V(f)$ is singular at the point $(0:0:1)$ no matter the other coefficients of $f$, and correspondingly $\operatorname{Disc}(f)$ vanishes identically; the polynomial $\varphi$ vanishes at those values of the other coefficients of $f$ for which the curve $V(f)$ is more singular than expected, i.e., has singularities worse than an ordinary double point at $(0:0:1)$. It seems that Salmon had a good idea of what he was talking about, but visibly it was so common to him that it did not require any kind of explanation. This knowledge however has then apparently been completely forgotten.

Salmon then uses formula (🦉) to derive various enumerative quantities for surfaces in $\mathbf{P}^3$ by elimination from the latter. In particular, he computes the number of bitangent planes passing through a fixed general point $p \in \mathbf{P}^3$. His method is to consider a pencil of planes passing through a fixed point $p'$ on the surface in question, chosen such that the tangent plane at $p'$ is a member of the pencil. This pencil contains a fixed number of planes tangent to the surface, among which the tangent plane at $p'$ counts with multiplicity 2 in general, and with greater multiplicity if it has some special feature, e.g., if it is a bitangent plane. In an appropriate setting, this multiplicity is the valuation in $T$ of the polynomial in (🦆), and the game is to understand the conditions that make it jump. It is maybe not so surprising that the techniques we use to follow this plan have the same flavour than those with which we obtain formula (🦆) in the first place: it is all based on a thorough study of the various homogeneity properties of the resultant (and as a special case, of the discriminant) and their interplay. We group these techniques under the concept of *reduced elimination*.

There are other well-known ways to compute the number of bitangent planes to a surface in $\mathbf{P}^3$ (see, e.g., [6]). Let us list a few points in favour of Salmon's technique. First of all it is the natural approach, and does not involve any trick. A concrete manifestation of this is that Salmon's technique gives more than a mere degree: it shows the existence of a node-couple hypersurface, the intersection of which with the surface under consideration is the locus of tangency points of bitangent planes (see Theorem (3.17)). Moreover, it works over an arbitrary algebraically closed field of characteristic zero, and can be carried out for a hypersurface in a projective space of arbitrary dimension, as we observe in the present text.

Indeed, we prove the following. Let $X \subset \mathbf{P}^n$ be a hypersurface. In a suitable system of homogeneous coordinates, it is defined by a homogeneous polynomial of the form

$$f(x_0, x_1, \ldots, x_n) = T x_0^{d-1} x_n + \sum_{k \geqslant 2} x_0^{d-k} f_k(x_1, \ldots, x_n),$$

and we prove in (2.17) that one has

$$\mathrm{Disc}(f) = T^2 \mathrm{Disc}(\bar{f}_2) \mathrm{Disc}(f_2)^2 \, \mathrm{redDisc}(f) \bmod T^3,$$

where $\bar{f}_2(x_1, \ldots, x_{n-1}) = f_2(x_1, \ldots, x_{n-1}, 0)$, and $\mathrm{redDisc}(f)$ is the *reduced discriminant* of $f$ with respect to the truncation at the order $d-2$ in the variable $x_0$, see Definition (2.2). From this we are able to show that, under suitable smoothness and transversality assumptions, the number of bitangent hyperplanes to $X$ passing through $n-2$ general points in $\mathbf{P}^n$ is

$$\frac{1}{2} d(d-1)^{n-2}(d-2)\Big(d \cdot \frac{(d-1)^n - 1}{d-2} - 3(n+1)\Big)$$

(see Theorem (3.24)).

In addition, we have included in the text two more projective enumerative computations, as further applications of reduced elimination theory. Namely, we give the computations of the respective numbers of flex-tangent hyperplanes and bitangent lines to a surface in $\mathbf{P}^3$, again following Salmon's ideas (note that the former computation had already been carried out in detail and in arbitrary dimension in [2]).

We also take the occasion to give a synthetic account of the basic theory of resultants and discriminants, in a way which we believe could be useful to the early XXI$^{\mathrm{st}}$ century classical algebraic geometer.

The article is organized as follows. In Section 1, we review the theory of the resultant including its various homogeneity properties, and introduce the reduced resultant following Zariski. In Section 2 we discuss in the same fashion the ordinary discriminant and its reduced version; in subsection 2.3 we prove Salmon's formula (🦆) and its version in arbitrary dimension. Section 3 is devoted to the computation by elimination from the latter formula of the number of bitangent hyperplanes to a smooth hypersurface, with emphasis on the surface case. The final Section 4 contains the additional computations of the flecnodal degree and of the number of bitangent lines to a surface in $\mathbf{P}^3$.

This text is a slightly modified version of a chapter of the forthcoming proceedings of the Seminar on Degenerations and enumeration of curves on surfaces held in Roma "Tor Vergata" 2015–2017. The main modification consists in the inclusion of subsection 3.1 on polarity to make the text more self-contained.

Th.D. wishes to thank Ragni Piene and Israel Vainsencher for their interest and a crucial suggestion at a prehistorical stage of this work. L.B. is grateful to Alexandru Dimca for useful discussions on Milnor number and for indicating the reference [9].

# 1 − Reduced resultant

Suppose that $n + 1$ homogeneous polynomials $f_1, \ldots, f_{n+1}$ in the variables $x_0, x_1, \ldots, x_n$ are given. They define a collection of $n+1$ hypersurfaces in a projective space $\mathbf{P}^n$, the intersection of which is expected to be empty if they are sufficiently general. The emptiness of this intersection is indeed characterized by the non-vanishing of the resultant $\mathrm{Res}(f_1, \ldots, f_{n+1})$ of these polynomials (see §1.1). Thus, the resultant characterizes those collections of polynomials that have a common root. The purpose of the *reduced resultant* is similar: collections of polynomials $f_1, \ldots, f_{n+1}$ are still considered but with the additional property that they already have a common non-trivial root; then the reduced resultant will characterize those collections of polynomials having an additional extra root. It has been introduced by Zariski [18]; a more complete and modern treatment can be found in [14].

In what follows, we use the resultant of multivariate homogeneous polynomials as developed in [13] (see also [5, 7, 10]).

**(1.1) Notation.** Let $\mathbf{k}$ be a commutative ring, $n \geqslant 1$ an integer, $\mathbf{x} := (x_0, x_1, \ldots, x_n)$ a sequence of indeterminates. Given a multi-index $\boldsymbol{\alpha} := (\alpha_0, \alpha_1, \ldots, \alpha_n) \in \mathbf{N}^n$, we denote by $\mathbf{x}^{\boldsymbol{\alpha}}$ the monomial $x_0^{\alpha_0} x_1^{\alpha_1} \ldots x_n^{\alpha_n}$ and set $|\boldsymbol{\alpha}| = \sum_{i=0}^{n} \alpha_i$.

## 1.1 − Inertia forms and the ordinary multivariate resultant

We follow the beautiful presentation in [17, Chapter XI] (beware that this Chapter on Elimination Theory has disappeared in later editions).

**(1.2) Saturation of a homogeneous ideal.** We recall the following for the convenience of the reader; see, e.g., [11, Lecture 5] or [12, Exercise II.5.10] for more details. The *saturation* of a homogeneous ideal $I \subset \mathbf{k}[\mathbf{x}]$ is the homogeneous ideal

$$\begin{aligned}\bar{I} &= \left\{ f \in \mathbf{k}[\mathbf{x}] : \forall i = 0, \ldots, n, \ \exists N_i \text{ s.t. } x_i^{N_i} f \in I \right\} \\ &= I : (x_0, \ldots, x_n)^{\infty}.\end{aligned}$$

For sufficiently large $m$, the graded pieces $I_m$ and $\bar{I}_m$ are equal. Moreover, for $I, J$ two homogeneous ideals, the following three propositions are equivalent:
(i) $\bar{I} = \bar{J}$;
(ii) $I_m = J_m$ for sufficiently large $m$;
(iii) $I \cdot \mathbf{k}[\mathbf{x}, \frac{1}{x_i}] = J \cdot \mathbf{k}[\mathbf{x}, \frac{1}{x_i}]$ for all $i = 0, \ldots, n$.
In other words, a subscheme $X \subset \mathbf{P}_{\mathbf{k}}^n$ is defined (scheme-theoretically) by a homogeneous ideal $I \subset \mathbf{k}[\mathbf{x}]$ if and only if the saturation $\bar{I}$ equals the homogeneous ideal $I_X$ of $X$.

In particular, if $\mathbf{k}$ is a field, a subscheme $X \subset \mathbf{P}_{\mathbf{k}}^n$ defined by a homogeneous ideal $I$ is empty if and only if the degree 0 piece $\bar{I}_0$ is non-zero — this is the homogeneous nullstellensatz (see also [4]). On the other hand, the non-emptiness of $X$ is equivalent to the existence of a point in $X(\mathbf{k}')$ for some finite field extension $\mathbf{k}'$ of $\mathbf{k}$.

When **k** is an arbitrary commutative ring, the vanishing of $\bar{I}_0$ is equivalent to the scheme-theoretic image of the map $X \to \operatorname{Spec}(\mathbf{k})$ being equal to the whole $\operatorname{Spec}(\mathbf{k})$. Indeed, the subscheme of $\operatorname{Spec}(\mathbf{k})$ defined by $\bar{I}_0$ coincides as a set with the image of $X \to \operatorname{Spec}(\mathbf{k})$ — this is the proof that projective morphisms are closed —, and moreover $\bar{I}_0$ defines the scheme-theoretic image of $X \to \operatorname{Spec}(\mathbf{k})$ (see, e.g., [13, §1] for more details).

**(1.3)** Let $d_1, \ldots, d_r$ be positive integers. For all $j = 1, \ldots, r$, we consider the generic homogeneous degree $d_j$ polynomial (in the variables **x**)

$$f_j := \sum_{|\boldsymbol{\alpha}| = d_j} u_{j,\boldsymbol{\alpha}} \mathbf{x}^{\boldsymbol{\alpha}}.$$

We set $A_{\mathbf{Z}} := \mathbf{Z}[u_{j,\boldsymbol{\alpha}} : j = 1, \ldots, r, |\boldsymbol{\alpha}| = d_j]$, so that $f_j \in A_{\mathbf{Z}}[\mathbf{x}]$ for all $j = 1, \ldots, r$.

**(1.4) Definition.** *An* inertia form *for the polynomials $f_1, \ldots, f_r$ is an element $T \in A_{\mathbf{Z}}$ such for all $i = 0, \ldots, n$, there exists $N_i \in \mathbf{N}$ such that $x_i^{N_i} T \in (f_1, \ldots, f_r)$.*

In other words, the inertia forms for $f_1, \ldots, f_r$ are the homogeneous elements of degree 0 of the saturation of the ideal $(f_1, \ldots, f_r)$ in $A_{\mathbf{Z}}[\mathbf{x}]$. The inertia forms for $f_1, \ldots, f_r$ form a prime ideal (see, e.g., [3, §2.1]) that we denote by $\mathfrak{J}_{\mathbf{Z}}$.

**(1.5) Theorem.** *Suppose that **k** is a field, and let $a_{j,\boldsymbol{\alpha}} \in \mathbf{k}$ for $1 \leqslant j \leqslant r$ and $|\boldsymbol{\alpha}| \leqslant d_j$, and consider the polynomials $f_j(a_{j,\boldsymbol{\alpha}}) \in \mathbf{k}[\mathbf{x}]$ obtained by evaluating each variable $u_{j,\boldsymbol{\alpha}}$ in $a_{j,\boldsymbol{\alpha}} \in \mathbf{k}$. The two following propositions are equivalent:*
*(i) the ideal $(f_1(a_{1,\boldsymbol{\alpha}}), \ldots, f_r(a_{r,\boldsymbol{\alpha}}))$ defines a non-empty subscheme of $\mathbf{P}^n_{\mathbf{k}}$;*
*(ii) for all $T \in \mathfrak{J}_{\mathbf{Z}}$, for all $j = 1, \ldots, r$ and $|\boldsymbol{\alpha}| \leqslant d_j$: $T(a_{j,\boldsymbol{\alpha}}) = 0$.*

This tells us that a given specialization to a field of the polynomials $f_j$ defines a non-empty subscheme if and only if all the constants in the saturation of $(f_1, \ldots, f_r)$ specialize to 0 in this specialization (see also [4]).

We emphasize that in general the subscheme of $\operatorname{Spec}(\mathbf{k})$ defined by the specialization of $\mathfrak{J}_{\mathbf{Z}}$ only coincides set-theoretically with the scheme-theoretic image of $X \to \operatorname{Spec}(\mathbf{k})$ (see [7, §3, Remarque 1] and [13, §1]), which is the reason why we assume that **k** is a field in Theorem (1.5). If **k** is an arbitrary commutative ring, what is indeed true is that the subscheme defined by the $f_j(a_{j,\alpha})$'s surjects onto $\operatorname{Spec}(\mathbf{k})$ as a set if and only if $\mathfrak{J}_{\mathbf{Z}} \otimes_{A_{\mathbf{Z}}} \mathbf{k}$ is contained in the nilradical $\sqrt{(0)}$ of **k**, but this says nothing more than Theorem (1.5).

**(1.6) Theorem.** *If $r = n+1$, the ideal $\mathfrak{J}_{\mathbf{Z}}$ is principal. Up to sign it has a single generator that we denote by $\operatorname{Res}_{d_1,\ldots,d_{n+1}} \in A_{\mathbf{Z}}$. The latter is an irreducible element of $A_{\mathbf{Z}}$. Moreover, for all $k \in [\![1, n+1]\!]$, $\operatorname{Res}_{d_1,\ldots,d_{n+1}}$ is homogeneous of degree $\prod_{j \neq k} d_j$ with respect to the coefficients of the polynomial $f_k$, i.e., with respect to the variables $u_{k,\boldsymbol{\alpha}}$, $|\boldsymbol{\alpha}| \leqslant d_k$ (all assumed to have weight equal to one).*

Let $\tilde{f}_1, \ldots, \tilde{f}_{n+1} \in \mathbf{k}[\mathbf{x}]$ be polynomials of respective degrees $d_1, \ldots, d_{n+1}$. They are specializations of $f_1, \ldots, f_{n+1} \in A_{\mathbf{Z}}[\mathbf{x}]$ for an appropriate canonical choice of $a_{j,\boldsymbol{\alpha}} \in \mathbf{k}$. We let $\operatorname{Res}(\tilde{f}_1, \ldots, \tilde{f}_{n+1}) \in \mathbf{k}$ (or $\operatorname{Res}_{d_1,\ldots,d_{n+1}}(\tilde{f}_1, \ldots, \tilde{f}_{n+1}) \in \mathbf{k}$, if we want to emphasize the dependency on the degrees) be the corresponding specialization of $\operatorname{Res}_{d_1,\ldots,d_{n+1}} \in A_{\mathbf{Z}}$. The multi-homogeneity property stated in the above theorem may then be rephrased as follows: for all $\lambda \in \mathbf{k}$,

$$\operatorname{Res}_{d_1,\ldots,d_{n+1}}(\tilde{f}_1, \ldots, \lambda \tilde{f}_k, \ldots, \tilde{f}_{n+1}) = \lambda^{d_1 \cdots d_{k-1} d_{k+1} \cdots d_{n+1}} \operatorname{Res}_{d_1,\ldots,d_{n+1}}(\tilde{f}_1, \ldots, \tilde{f}_k, \ldots, \tilde{f}_{n+1}).$$

The sign indeterminacy in the definition of $\operatorname{Res}_{d_1,\ldots,d_{n+1}}$ is usually removed by imposing the equality $\operatorname{Res}(x_0^{d_0},\ldots,x_n^{d_n}) = 1$.

**(1.7) Divisibility property** (see, e.g., [13, §5.6]). Let $\tilde{f}_1,\ldots,\tilde{f}_{n+1}$ and $\tilde{g}_1,\ldots,\tilde{g}_{n+1}$ be two sequences of homogeneous polynomials in $\mathbf{k}[\mathbf{x}]$ such that we have the inclusion of ideals of $\mathbf{k}[\mathbf{x}]$

$$(\tilde{g}_1,\ldots,\tilde{g}_{n+1}) \subset (\tilde{f}_1,\ldots,\tilde{f}_{n+1}).$$

Then $\operatorname{Res}(\tilde{f}_1,\ldots,\tilde{f}_{n+1})$ divides $\operatorname{Res}(\tilde{g}_1,\ldots,\tilde{g}_{n+1})$ in $\mathbf{k}$.

**(1.8) Weight properties.** Besides its ordinary multi-homogeneity property (see Theorem (1.6)) the resultant has other homogeneous structures that we call "weight properties" to emphasize that the grading of the coefficient ring $A_\mathbf{Z}$ is not the standard one.

With the notation (1.3), let $k$ be an integer in $[\![0,n]\!]$ and define a new grading on $A_\mathbf{Z} = \mathbf{Z}[u_{j,\boldsymbol{\alpha}} : j = 1,\ldots,r, |\boldsymbol{\alpha}| = d_j]$ by setting

$$(1.8.1) \qquad \operatorname{weight}(u_{j,\boldsymbol{\alpha}}) = \alpha_k.$$

Then, $\operatorname{Res}_{d_1,\ldots,d_{n+1}}$ is homogeneous of degree $d_1 d_2 \cdots d_{n+1}$ with respect to this new grading (see [13, §5.13.2]). The Bezout theorem, which counts the number of roots of a finite complete intersection scheme in a projective space and thus is the mother of all statements in projective enumerative geometry, can be deduced from this property.

Another interesting weight property of the resultant is obtained by grading $A_\mathbf{Z}$ with

$$(1.8.2) \qquad \operatorname{weight}(u_{j,\boldsymbol{\alpha}}) = d_j - \alpha_k.$$

In this case, $\operatorname{Res}_{d_1,\ldots,d_{n+1}}$ is homogeneous of degree $n d_1 d_2 \cdots d_{n+1}$ (this follows from [13, §5.13] as well).

**(1.9) Further weight properties.** One may combine the previous weight properties of (1.8) with the standard homogeneity of the resultant in Theorem (1.6) to obtain further weight properties. The argument is as follows.

Assume that the resultant is homogeneous of degree $\delta$ for some grading on $A_\mathbf{Z}$, and let $\operatorname{weight}(u_{j,\boldsymbol{\alpha}}) = w_{j,\boldsymbol{\alpha}}$. Since the resultant is a homogeneous polynomial of degree $\prod_{j \neq k} d_j$ with respects to the variables $u_{k,\boldsymbol{\alpha}}$ (for the standard grading), a shift by $r$ in the weights of all the variables $u_{k,\boldsymbol{\alpha}}$ for some $k$ induces a shift by $r \cdot \prod_{j \neq k} d_j$ in the degree of the resultant.

Let $r_1,\ldots,r_{n+1} \in \mathbf{Z}$, and consider the new grading on $A_\mathbf{Z}$ defined by setting $\operatorname{weight}(u_{j,\boldsymbol{\alpha}}) = w_{j,\boldsymbol{\alpha}} + r_j$. The resultant is homogeneous with respect to this new grading, of degree

$$\delta + \sum_{1 \leqslant k \leqslant n+1} \left( r_k \cdot \prod_{j \neq k} d_j \right).$$

## 1.2 – The reduced resultant

We shall now explain how to adapt the ideas of the previous paragraph to develop the theory of the reduced resultant. We refer to [18] and [14] for the details and proofs. Somehow, this is a generalization of the following toy example.

**(1.10) Example** (projection of a complete intersection from one of its points). Let $f, g \in \mathbf{k}[\mathbf{x}]$ be two homogeneous polynomials of degrees $a$ and $b$, defining a complete intersection $X \subset \mathbf{P}^n$, and suppose one wants to project $X$ from a point $p_0 \in \mathbf{P}^n$. Assume for simplicity that $\mathbf{k}$ is an algebraically closed field. We may take $p_0 = (1 : 0 : \ldots : 0)$. Then one considers the two polynomials

$$
\begin{aligned}
f(t, x_1, \ldots, x_n) &= f_0 t^a + f_1 t^{a-1} + \cdots + f_a \\
\text{and} \quad g(t, x_1, \ldots, x_n) &= g_0 t^b + g_1 t^{b-1} + \cdots + g_b
\end{aligned}
\tag{1.10.1}
$$

in $\mathbf{k}[x_1, \ldots, x_n][t]$ for all $(x_1 : \ldots : x_n) \in \mathbf{P}^{n-1}$. (We are abusing notation here, as one should consider instead the two polynomials $f(t, sx_1, \ldots, sx_n)$ and $g(t, sx_1, \ldots, sx_n)$ that are homogeneous in the couple of variables $(s, t)$).

If $p_0 \notin X$, the point $(x_1 : \ldots : x_n) \in \mathbf{P}^{n-1}$ belongs to the projection of $X$ from $p_0$ if and only if the two polynomials in (1.10.1) have a common root in $\mathbf{P}^1$, hence the equation of the projection is given by

$$\operatorname{Res}_{a,b}(f, g) \in \mathbf{k}[x_1, \ldots, x_n],$$

which is homogeneous of degree $ab$ in the variables $(x_1, \ldots, x_n)$ by (1.8), as the polynomials $f_i$ are homogeneous of degree $i$ in $\mathbf{k}[x_1, \ldots, x_n]$.

On the other hand, if $p_0 \in X$ then, letting $a'$ and $b'$ be the respective multiplicities of $p_0$ in the hypersurfaces $V(f)$ and $V(g)$, one has

$$f_0 = \cdots = f_{a'-1} = g_0 = \cdots = g_{b'-1} = 0,$$

so that (1.10.1) becomes

$$
\begin{aligned}
f(t, x_1, \ldots, x_n) &= f_{a'} t^{a-a'} + \cdots + f_a = {}^{\flat}\!f \\
\text{and} \quad g(t, x_1, \ldots, x_n) &= g_{b'} t^{b-b'} + \cdots + g_b = {}^{\flat}\!g.
\end{aligned}
\tag{1.10.2}
$$

It follows that the equation of the projection of $X$ from $p_0$ is given by

$$\operatorname{Res}_{a-a', b-b'}({}^{\flat}\!f, {}^{\flat}\!g) \in \mathbf{k}[x_1, \ldots, x_n].$$

We shall see later on that this polynomial is the *reduced resultant* of $f$ and $g$ truncated at orders $a - a'$ and $b - b'$ respectively, as polynomials in the variable $t$. It is a homogeneous polynomial of degree $ab - a'b'$ in the variables $(x_1, \ldots, x_n)$: this may be seen using (1.8) and (1.9). Indeed, the coefficient of ${}^{\flat}\!f$ (resp. ${}^{\flat}\!g$) in $t^i$ is a homogeneous polynomial in $(x_1, \ldots, x_n)$ of degree $a - i = \deg_t({}^{\flat}\!f) - i + a'$ (resp. $b - i = \deg_t({}^{\flat}\!g) - i + b'$). Therefore, the argument of (1.9) applied to (1.8.2) gives that $\operatorname{Res}_{a-a',b-b'}({}^{\flat}\!f, {}^{\flat}\!g)$ is homogeneous of degree

$$(a - a')(b - b') + a'b + b'a = ab - a'b'$$

as we had announced. This weight property is a particular case of (1.16) which applies to reduced resultants in general.

**(1.11)** Let $d_1, \ldots, d_{n+1}$ be positive integers, and consider for all $j \in [\![1, n+1]\!]$ the generic homogeneous degree $d_j$ polynomial

$$f_j = \sum_{|\boldsymbol{\alpha}| = d_j} u_{j, \boldsymbol{\alpha}} \mathbf{x}^{\boldsymbol{\alpha}} = \sum_{k=0}^{d_j} x_0^{d_j - k} f_{j,k}(x_1, \ldots, x_n) \in A_{\mathbf{Z}}[\mathbf{x}].$$

So, $f_{j,k}$ is a homogeneous polynomial in $A_{\mathbf{Z}}[x_1, \ldots, x_n]$ of degree $k$. Fix an integer $s_j \in [\![1, d_j]\!]$ for all $j$. The *truncation* of $f_j$ at order $d_j - s_j$ with respect to $x_0$ is the polynomial

$$h_j = \sum_{k=s_j}^{d_j} x_0^{d_j - k} f_{j,k} = x_0^{d_j - s_j} f_{j,s_j} + \cdots + x_0 f_{j, d_j - 1} + f_{j, d_j} \in A_{\mathbf{Z}}[\mathbf{x}].$$

The purpose of reduced elimination theory is the study of inertia forms of the truncations of polynomials at some given orders; essentially, this can be done with the same strategy as in the classical case we recalled in §1.1.

**(1.12) Theorem** ([18, Theorem 6 and §8] and [14, Theorem II.0.5 and §IV.0]). *Assume that $d_j > s_j$ for some $j \in [\![1, n+1]\!]$. The ideal of reduced inertia forms*

$$\mathcal{Q}_{\mathbf{Z}} = \big((h_1, \ldots, h_{n+1}) : (x_1, \ldots, x_n)^\infty\big) \cap A_{\mathbf{Z}}$$

*is a prime and principal ideal of $A_{\mathbf{Z}}$. The* reduced resultant, *denoted*

$$\mathrm{redRes}_{d_1, \ldots, d_{n+1}}^{s_1, \ldots, s_{n+1}} \in A_{\mathbf{Z}},$$

*is defined, up to sign, as the generator of $\mathcal{Q}_{\mathbf{Z}}$; it is therefore an irreducible element of $A_{\mathbf{Z}}$.*

*Moreover, if $d_j > s_j$ for at least two distinct integers $j, j' \in [\![1, n+1]\!]$, then for all $i \in [\![1, n+1]\!]$ the reduced resultant is a homogeneous polynomial of degree*

$$\frac{d_1 d_2 \cdots d_{n+1}}{d_i} - \frac{s_1 s_2 \cdots s_{n+1}}{s_i}$$

*with respect to the coefficients of the polynomial $h_i$, i.e., with respect to the coefficients $u_{i, \boldsymbol{\alpha}}$ such that $|\boldsymbol{\alpha}| = d_i$ and $\alpha_0 \leqslant d_i - s_i$.*

*If there is only one integer $j \in [\![1, n+1]\!]$ such that $d_j > s_j$, then the reduced resultant is equal to the resultant of the polynomials $h_1, \ldots, h_{j-1}, h_{j+1}, \ldots, h_{n+1}$.*

Note that in the above statement the ideal $(h_1, \ldots, h_{n+1})$ is saturated with respect to $(x_1, \ldots, x_n)$ — which is the defining ideal of the point $(1 : 0 : \ldots : 0)$ — instead of $(x_0, x_1, \ldots, x_n)$, and that the polynomials $h_1, \ldots, h_{n+1}$ are not homogeneous in the set of variables $(x_1, \ldots, x_n)$.

Of course the reduced resultant only depends on the coefficients of the generic truncated polynomials $h_1, \ldots, h_{n+1}$, and not on all the coefficients of the polynomial $f_1, \ldots, f_{n+1}$. We will often denote it by $\mathrm{redRes}(h_1, \ldots, h_{n+1})$ without printing the integers $d_i$ and $s_i$, that are implicitly given by the polynomials $h_1, \ldots, h_{n+1}$. The reduced resultant of a collection of polynomials $\tilde{h}_1, \ldots, \tilde{h}_{n+1} \in \mathbf{k}[\mathbf{x}]$ is defined as the corresponding specialization of the generic reduced resultant; it is an element in $\mathbf{k}$ denoted by $\mathrm{redRes}(\tilde{h}_1, \ldots, \tilde{h}_{n+1})$. The sign indeterminacy in the definition of the reduced resultant can be removed by means of Theorem (1.15), once the sign of the multivariate resultant is chosen.

**(1.13) Vanishing of the reduced resultant.** The reduced resultant $\mathrm{redRes}(h_1, \ldots, h_{n+1})$ is a polynomial in the ring of coefficients of the polynomials $h_j$, $j = 1, \ldots, n+1$, i.e., in the ring

$$A_{\mathbf{Z}}[u_{i, \boldsymbol{\alpha}} : i = 1, \ldots, n+1, |\boldsymbol{\alpha}| = d_i, \alpha_0 \leqslant d_i - s_i].$$

Its vanishing on an algebraically closed field $\mathbf{k}$ characterizes those collections of hypersurfaces $\tilde{h}_1, \ldots, \tilde{h}_{n+1}$ of $\mathbf{P}_{\mathbf{k}}^n$ that have a further intersection point, infinitely near or not, besides the origin $(1 : 0 : \ldots : 0)$, i.e., those collections that satisfy to one of the two following conditions:

(a) the hypersurfaces $\tilde{h}_1, \ldots, \tilde{h}_{n+1}$ intersect at an additional point which is not $(1 : 0 : \ldots : 0)$;

(b) the polynomials $\tilde{f}_{j,s_j}$, $j = 1, \ldots, n+1$, have a common root in $\mathbf{P}_{\mathbf{k}}^{n-1}$, which means that the hypersurfaces $\tilde{h}_1, \ldots, \tilde{h}_r$ have a common principal tangent at the point $(1 : 0 : \ldots : 0)$.

These properties are proved in [18, Theorem 3.1, 3.2 and 3.3], and in [14, Proposition I.1].

## 1.3 – Generalized weight properties

In [18] Zariski showed that the reduced resultant can be computed from its corresponding resultant. To obtain this property, he introduced a generalization of the grading (1.8.1) and, although the resultant is no longer homogeneous with respect to this new grading, he proved that its graded part of smallest degree is connected to the reduced resultant.

We maintain the notation of §1.2.

**(1.14) The Zariski grading.** We define a grading on $A_{\mathbf{Z}} = \mathbf{Z}[u_{j,\boldsymbol{\alpha}}]$ by assigning for all $j$

$$\text{weight}(u_{j,\boldsymbol{\alpha}}) = \begin{cases} 0 & \text{if } \alpha_0 < d_j - s_j \\ \alpha_0 - d_j + s_j & \text{otherwise,} \end{cases}$$

and weight 0 to the constants. We find it helpful to visualize this definition as follows:

$$f_j = \underbrace{x_0^{d_j} f_{j,0}}_{\text{coeffs have weight}=s_j} + \cdots + \underbrace{x_0^{d_j-s_j+1} f_{j,s_j-1}}_{\text{coeffs have weight}=1} + \underbrace{x_0^{d_j-s_j} f_{j,s_j} + \cdots + f_{j,d_j}}_{\text{coeffs have weight}=0}.$$

Note in particular that the coefficients $u_{j,\boldsymbol{\alpha}}$ whose weight is equal to 0 in this grading are exactly the coefficients of the truncation $h_j$ of the polynomial $f_j$, respectively. The grading (1.8.1) introduced in (1.8) is a particular case of a Zariski grading (corresponding to $s_j = d_j$ for all $j$), which explains the terminology "generalized weight properties".

The main property of the Zariski grading is that it allows the computation of the reduced resultant of $h_1, \ldots, h_{n+1}$ (the truncations of $f_1, \ldots, f_{n+1}$ at the orders $d_1 - s_1, \ldots, d_{n+1} - s_{n+1}$, respectively) from the resultant of $f_1, \ldots, f_{n+1}$. To this end, we introduce one more notation: for all $j = 1, \ldots, n+1$, we let $g_j$ be the quotient of the Euclidean division of $f_j$ by $x_0^{d_j-s_j}$ in $A_{\mathbf{Z}}[x_1, \ldots, x_n][x_0]$, i.e.,

$$g_j = \frac{1}{x_0^{d_j-s_j}} \sum_{k \leqslant s_j} x_0^{d_j-k} f_{j,k} = \sum_{k=0}^{s_j} x_0^{s_j-k} f_{j,k} = x_0^{s_j} f_{j,0} + \cdots + x_0 f_{j,s_j-1} + f_{j,s_j}.$$

The polynomial $g_j$ is a generic degree $s_j$ homogeneous polynomial in the set of variables $\mathbf{x}$. Beware that $f_j$ does not equal $h_j + x_0^{d_j-s_j} g_j$, as $f_{j,s_j}$ appears in both $h_j$ and $g_j$.

**(1.15) Theorem** ([18, Theorem 5.1 and Theorem 5.2] and [14, Lemme IV.1.6]). *The nonzero homogeneous piece of lowest degree with respect to the Zariski grading of $\mathrm{Res}(f_1, \ldots, f_{n+1})$ is of degree $s_1 s_2 \cdots s_{n+1}$. Moreover, it equals the product*

$$\mathrm{Res}(g_1, \ldots, g_{n+1}) \cdot \mathrm{redRes}(h_1, \ldots, h_{n+1})$$

*if $s_j < d_j$ for at least two distinct integers $j, j' \in [\![1, n+1]\!]$, and the product*

$$\mathrm{Res}(g_1, \ldots, g_{n+1}) \cdot \mathrm{redRes}(h_1, \ldots, h_{n+1})^{d_j-s_j}$$

*if $s_j < d_j$ and $s_i = d_i$ for all integers $i \in [\![1, n+1]\!] \setminus \{j\}$.*

We notice that since the coefficients of $h_1, \ldots, h_{n+1}$ all have weight 0, the valuation of $\mathrm{Res}(f_1, \ldots, f_{n+1})$ is the degree of $\mathrm{Res}(g_1, \ldots, g_{n+1})$, which is equal to $\prod_{1 \leqslant j \leqslant n+1} s_j$ by (1.8).

**(1.16)** Suppose that the ring $A_{\mathbf{Z}}$ is graded by (1.8.2), with $k = 0$, which means that the coefficients of $f_{j,l}$ all have weight $l$. Then we know by (1.8) that the resultant of $f_1, \ldots, f_{n+1}$ is homogeneous of degree $nd_1 \cdots d_{n+1}$. And by the same argument, the resultant of $g_1, \ldots, g_{n+1}$ is homogeneous of degree $ns_1 \cdots s_{n+1}$. Therefore, we deduce from Theorem (1.15) that the reduced resultant of $h_1, \ldots, h_{n+1}$ is homogeneous of degree $n(d_1 \cdots d_{n+1} - s_1 \cdots s_{n+1})$ with respect to the grading (1.8.2), providing there are at least two distinct integers $j$ such that $s_j < d_j$. If $s_j < d_j$ and $s_i = d_i$ for all integers $i \in [\![1, n+1]\!] \setminus \{j\}$, then the reduced resultant is homogeneous of degree $nd_1 \ldots d_{n+1}/d_j$. We notice that this homogeneity property of the reduced resultant applies to give the degree of $\mathrm{Res}(^\flat f, ^\flat g)$ in Example (1.10).

## 2 – Reduced discriminant and Salmon formula

In this Section we give a rigorous proof to formula (). This is done by introducing the concept of reduced discriminant. We begin with a quick recap on the ordinary discriminant of a hypersurface, following [3, § 4]; see also [7] and [10, Chapter 13, §D].

### 2.1 – Discriminant of a homogeneous polynomial

**(2.1)** Let $d$ be a positive integer, and consider the generic homogeneous degree $d$ polynomial $f = \sum_{|\boldsymbol{\alpha}|=d} u_{\boldsymbol{\alpha}} \mathbf{x}^{\boldsymbol{\alpha}}$ in $n+1$ variables $\mathbf{x} = (x_0, \ldots, x_n)$. We set $A_{\mathbf{Z}} = \mathbf{Z}[u_{\boldsymbol{\alpha}} : |\boldsymbol{\alpha}| = d]$. For all $i = 0, \ldots, n$ we let $\partial_i$ denote derivation with respect to the variable $x_i$.

**(2.2) Definition.** *There is a unique element $\mathrm{Disc}_d(f) \in A_{\mathbf{Z}}$ (often simply denoted by $\mathrm{Disc}(f)$) such that*

$$(2.2.1) \qquad d^{a(n,d)} \mathrm{Disc}_d(f) = \mathrm{Res}(\partial_0 f, \ldots, \partial_n f)$$

*in $A_{\mathbf{Z}}$, where $a(n,d) = \frac{(d-1)^{n+1} - (-1)^{n+1}}{d} \in \mathbf{Z}$. It is homogeneous of degree $(n+1)(d-1)^n$ with respect to the coefficients of the polynomial $f$, i.e., with respect to the indeterminates $u_{\boldsymbol{\alpha}}$, $|\boldsymbol{\alpha}| = d$.*

*For a homogeneous degree $d$ polynomial $\tilde{f} \in \mathbf{k}[\mathbf{x}]$, we define the* discriminant $\mathrm{Disc}(\tilde{f}) \in \mathbf{k}$ *of $\tilde{f}$ as the specialization $\sigma(\mathrm{Disc}(f)) \in \mathbf{k}$, where $\sigma : A_{\mathbf{Z}} \to \mathbf{k}$ is the unique specialization morphism mapping $f$ to $\tilde{f}$.*

**(2.3) Proposition.** *The ideal of inertia forms*

$$(\partial_0 f, \partial_1 f, \ldots, \partial_n f, f) : (x_0, \ldots, x_n)^\infty$$

*is a prime and principal ideal in $A_{\mathbf{Z}}$. It is generated by the discriminant $\mathrm{Disc}(f)$, which is therefore an irreducible polynomial in $A_{\mathbf{Z}}$.*

**(2.4) Theorem** (smoothness criterion)**.** *Suppose $\mathbf{k}$ is an algebraically closed field, and consider a degree $d$ homogeneous polynomial $f \in \mathbf{k}[\mathbf{x}]$. The following are equivalent:*
*(i) the hypersurface $V(f) \subset \mathbf{P}_{\mathbf{k}}^n$ is smooth;*
*(ii) $\mathrm{Disc}(f) \neq 0$.*

**(2.5)** Denote by $\bar{f}$ the polynomial $f(0, x_1, \ldots, x_n) \in A_{\mathbf{Z}}[x_1, \ldots, x_n]$; it is the equation of the hypersurface in $\mathbf{P}^{n-1}$ cut out by $V(f)$ on the hyperplane $V(x_0) \subset \mathbf{P}^n$, of which we think as the hyperplane at infinity. We have the following identity in $A_{\mathbf{Z}}$, somehow reminiscent of the Euler identity:

$$(2.5.1) \qquad d^{(d-1)^n} \mathrm{Res}(\partial_1 f, \ldots, \partial_n f, f) = \mathrm{Res}(\partial_0 f, \ldots, \partial_n f) \cdot \mathrm{Res}(\partial_1 \bar{f}, \ldots, \partial_n \bar{f}).$$

Note that at the right-hand-side of this identity, the first (resp. second) factor is the resultant of $n+1$ (resp. $n$) polynomials in $n+1$ (resp. $n$) variables. Since $a(n,d) + a(n-1,d) = (d-1)^n$, (2.5.1) is equivalent to

$$(2.5.2) \qquad \mathrm{Res}(\partial_1 f, \ldots, \partial_n f, f) = \mathrm{Disc}(f) \cdot \mathrm{Disc}(\bar{f}).$$

Let $\mathbf{k}$ be an algebraically closed field, and consider a degree $d$ homogeneous polynomial $f \in \mathbf{k}[\mathbf{x}]$. By Theorem (2.4), the vanishing of $\mathrm{Disc}(\bar{f})$ is equivalent to the hyperplane section at infinity $V(f) \cap V(x_0)$ being singular. For a general $f$ such that $\mathrm{Disc}(\bar{f}) = 0$, the hypersurface $V(f)$ is non-singular and tangent to the hyperplane $V(x_0)$.

**(2.6)** Similarly to the resultant, the discriminant is also homogeneous under the specific grading of the coefficient ring $A_{\mathbf{Z}}$ introduced in (1.8). More precisely, let $k$ be an integer in $[\![0, n]\!]$.

If $A_{\mathbf{Z}}$ is graded with the rule $\mathrm{weight}(u_\alpha) = \alpha_k$, then the discriminant of $f$ is homogeneous of degree $d(d-1)^n$. This follows in a straightforward manner from (2.5.2) and the corresponding weight property of the resultant. One may use this result to compute the degree of the dual to a smooth hypersurface in $\mathbf{P}^{n+1}$, using the approach of Section 3, see (3.16).

If $A_{\mathbf{Z}}$ is graded with the rule $\mathrm{weight}(u_\alpha) = d - \alpha_k$ then the discriminant of $f$ is homogeneous of degree $nd(d-1)^n$. This weight property is easily deduced from the invariance of the discriminant under linear change of coordinates for which we refer the reader to [3, Proposition 4.13].

In turn, one may reproduce the argument given in (1.9) to deduce further weight properties from the two latter results and the standard homogeneity property of the discriminant stated in Definition (2.2).

## 2.2 – The reduced discriminant

**(2.7)** We write the generic homogeneous degree $d$ polynomial as

$$f = x_0^d f_0 + \cdots + x_0^{d-s} f_s + \cdots + f_d$$

where each $f_k$ is homogeneous of degree $k$ in the variables $x_1, \ldots, x_n$. We choose an integer $s \in [\![2, d-1]\!]$ and consider the truncation of $f$ at order $d-s$ with respect to $x_0$; we set

$$h = x_0^{d-s} f_s + x_0^{d-s-1} f_{s+1} + \cdots + f_d,$$
$$g = x_0^s f_0 + x_0^{s-1} f_1 + \cdots + f_s.$$

The polynomial $h$ is of degree $d$ and of valuation $s$ with respect to the variables $x_1, \ldots, x_n$, and its partial derivatives with respect to the variables $x_1, \ldots, x_n$ are all of degree $d-1$ and of valuation $s-1$, so that the reduced resultant

$$\mathrm{redRes}(\partial_1 h, \ldots, \partial_n h, h) = \mathrm{redRes}_{d-1, \ldots, d-1, d}^{s-1, \ldots, s-1, s}(\partial_1 h, \ldots, \partial_n h, h)$$

is well defined.

**(2.8) Proposition.** *With the above notation, the discriminants* $\mathrm{Disc}(f_s)$ *and* $\mathrm{Disc}(f_d)$ *both divide the reduced resultant* $\mathrm{redRes}(\partial_1 h, \ldots, \partial_n h, h)$.

*Proof.* Both discriminants $\mathrm{Disc}(f_s)$ and $\mathrm{Disc}(f_d)$ are irreducible as elements of $A_{\mathbf{Z}}$. In addition, their vanishing implies the vanishing of the reduced resultant by (1.13). Indeed, the vanishing of $\mathrm{Disc}(f_d)$ implies the existence of a common root at infinity ($x_0 = 0$) of the polynomial system $\partial_1 h = \cdots = \partial_n h = h = 0$. In the same way, the vanishing of $\mathrm{Disc}(f_s)$ implies the existence of a common root of the polynomial system $\partial_1 h = \cdots = \partial_n h = h = 0$ infinitely near to the point $(1 : 0 : \ldots : 0)$. We thus conclude that $\mathrm{Disc}(f_s)$ and $\mathrm{Disc}(f_d)$ both divide $\mathrm{redRes}(\partial_1 h, \ldots, \partial_n h, h)$, which ends the proof.

Alternatively, this proposition can be proved by means of inertia forms, as follows. By Theorem (1.12) the reduced resultant $\mathrm{redRes}(\partial_1 h, \ldots, \partial_n h, h)$ belongs to the ideal of inertia forms $(\partial_1 h, \ldots, \partial_n h, h) : (x_1, \ldots, x_n)^\infty$. Therefore, for all integer $i \in [\![1, n]\!]$ there exists an integer $N_i$ such that

(2.8.1) $\qquad\qquad x_i^{N_i} \mathrm{redRes}(\partial_1 h, \ldots, \partial_n h, h) \in (\partial_1 h, \ldots, \partial_n h, h).$

Specializing the variable $x_0$ to 0 in (2.8.1), we immediately get that

$$x_i^{N_i} \mathrm{redRes}(\partial_1 h, \ldots, \partial_n h, h) \in (\partial_1 f_d, \ldots, \partial_n f_d, f_d)$$

for all $i \in [\![1, n]\!]$, from which we deduce that $\mathrm{redRes}(\partial_1 h, \ldots, \partial_n h, h)$ belongs to the ideal of inertia forms $(\partial_1 f_d, \ldots, \partial_n f_d, f_d) : (x_1, \ldots, x_n)^\infty$. By Proposition (2.3), this ideal is generated by the discriminant of the polynomial $f_d$, so $\mathrm{Disc}(f_d)$ divides $\mathrm{redRes}(\partial_1 h, \ldots, \partial_n h, h)$. A similar argument, but slightly more technical, can be used to show that $\mathrm{Disc}(f_s)$ divides $\mathrm{redRes}(\partial_1 h, \ldots, \partial_n h, h)$: see [14, Lemme I.1.3]. □

In the notation of (2.5), $f_d = \bar{f}$ and $f_s = \bar{g}$. Observe that $f_d$ and $f_s$ are generic homogeneous polynomials of degree $d$ and $k$ respectively, in the variables $(x_1, \ldots, x_n)$. Proposition (2.8) leads to the following definition.

**(2.9) Definition.** *The* reduced discriminant *of $f$ with respect to the truncation at order $d - s$ with respect to the variable $x_0$, denoted* $\mathrm{redDisc}_d^s(h)$, *or simply* $\mathrm{redDisc}(h)$, *is defined by the equality*

(2.9.1) $\qquad\qquad \mathrm{redRes}(\partial_1 h, \ldots, \partial_n h, h) = \mathrm{Disc}(f_d) \, \mathrm{Disc}(f_s) \, \mathrm{redDisc}(h) \in A_{\mathbf{Z}}.$

The identity (2.9.1) should be compared to (2.5.2). Beware that the reduced discriminant is not merely the reduced resultant of all the partial derivatives, because of the division by the factor $\mathrm{Disc}(f_s)$ (see (2.10) for comments about this factor).

A first consequence of this definition is that the reduced discriminant is a homogeneous polynomial in $A_{\mathbf{Z}}$, for the standard grading $\mathrm{weight}(u_\alpha) = 1$ for all $\alpha$, of degree

(2.9.2) $\qquad\qquad (n+1)\big[(d-1)^n - (s-1)^n\big] - 2n(s-1)^{n-1}.$

This follows from a straightforward computation since the degree of the other quantities in (2.9.1) are known.

We also note that $\mathrm{redDisc}(h)$ is a primitive polynomial in $A_{\mathbf{Z}}$ (i.e., the greatest common divisor of its coefficients equals 1), because $\mathrm{redRes}(\partial_1 h, \ldots, \partial_n h, h)$ is a primitive polynomial by [3, Proposition 4.24].

**(2.10) Vanishing of the reduced discriminant.** By definition, the generic hypersurface $V(h) \subset \mathbf{P}^n$ has an ordinary $s$-fold point at the origin $(1:0:\ldots:0)$. By (1.13) and (2.9), the vanishing of the product $\mathrm{Disc}(f_d) \cdot \mathrm{Disc}(f_s) \cdot \mathrm{redDisc}(h)$ corresponds to one of the two following properties:

(a) the existence of a common root to the polynomials $\partial_1 h, \ldots, \partial_n h, h$ which is not the origin,

(b) the existence of a common root to the polynomials $\partial_1 f_s, \ldots, \partial_n f_s, f_s$, equivalently a common root to the polynomials $\partial_1 f_s, \ldots, \partial_n f_s$ by Euler identity (assuming $d \neq 0$, which we do).

As Proposition (2.8) tells us, in each of these two cases there is a codimension one component in the space of hypersurfaces $V(h) \subset \mathbf{P}^n$ that can be factored out, namely the zero locus of $\mathrm{Disc}(f_d)$ in case (a), and that of $\mathrm{Disc}(f_s)$ in case (b). The factor $\mathrm{Disc}(f_d)$ is somehow artificial, corresponding as in (2.5) to our definition of the reduced discriminant with $\mathrm{redRes}(\partial_1 h, \ldots, \partial_n h, h)$ and not $\mathrm{redRes}(\partial_0 h, \partial_1 h, \ldots, \partial_n h)$. The reason why we proceed this way is that it lets us define the reduced discriminant as a primitive polynomial with integer coefficients without dealing with possible constant factors similar to $d^{a(n,d)}$ in (2.2.1). The factor $\mathrm{Disc}(f_s)$ on the other hand is indeed meaningful, as it characterizes those $h$ for which the tangent cone to $V(h)$ at the origin is a cone over a singular degree $s$ hypersurface (with the convention that a hypersurface of degree $s' > s$ is a singular degree $s$ hypersurface).

There is an interesting connection to Milnor number. Assume that $h$ defines a hypersurface with isolated singularities, so that its total Milnor number $\mu(h)$, and hence its Milnor number at the origin $\mu_0(h)$ are well defined. For a general such polynomial $h$, we have $\mu(h) = \mu_0(f_s) = (s-1)^n$ and property (b) is equivalent to the condition $\mu_0(h) > (s-1)^n$; see [9, Theorem 1]. Therefore, it follows that if $\mathrm{Disc}(f_d) \neq 0$, then $\mu(h) > (s-1)^n$ if and only $\mathrm{Disc}(f_s) \, \mathrm{redDisc}(h) = 0$. Pushing further, we see that if $\mathrm{Disc}(f_d) \neq 0$ and $\mathrm{Disc}(f_s) \neq 0$, then $\mathrm{redDisc}(h)$ vanishes if and only if $\mu(h) > \mu_0(h) = (s-1)^n$.

To conclude, the locus of those $h$ such that $V(h)$ either has a singularity off the origin, or has a singularity at the origin worse than an ordinary $s$-fold point, has several irreducible components, all of codimension 1: one is defined by the vanishing of $\mathrm{Disc}(f_s)$, and the reduced discriminant defines the others. It is plausible that there is only one such other component, i.e., that the reduced discriminant is irreducible, see (2.15).

**(2.11) Example.** As an illustrative example, we consider the case $n = 1$ and set

$$f = \sum_{i=0}^{d} a_i x_0^{d-i} x_1^i = x_0^d f_0 + x_0^{d-1} f_1 + \cdots + x_0^{d-s} f_s + \cdots f_d,$$

so that $f_i = a_i x_1^i$ for all $i = 0, \ldots, d$. We consider

$$h = \sum_{i=s}^{d} a_i x_0^{d-i} x_1^i = \sum_{i=s}^{d} x_0^{d-i} f_i = x_1^s \left( a_s x_0^{d-s} + a_{s+1} x_0^{d-s-1} x_1 + \cdots + a_d x_1^{d-s} \right),$$

the truncation of $f$ at order $d - s > 0$ with respect to $x_0$. Then, setting $h = x_1^s \cdot {}^\flat h$ it is easy to check that $\mathrm{redDisc}(h) = \pm \mathrm{Disc}({}^\flat h)$, which is therefore an irreducible polynomial of degree $2(d - s - 1)$ in the coefficients of $h$ (compare with the degree formula (2.9.2) in this setting). Actually, using the weight properties in (2.6), we can deduce weight properties of $\mathrm{redDisc}(h)$.

To be more precise, suppose that $A_{\mathbf{Z}}$ is graded with the rule $\mathrm{weight}(a_i) = \max(0, i - s)$, then the reduced discriminant $\mathrm{redDisc}(h)$ is homogeneous of degree $(d - s)(d - s - 1)$ by (2.6). And similarly, if $A_{\mathbf{Z}}$ is graded with the rule $\mathrm{weight}(a_i) = d - i$ then the same conclusion holds.

We can generalize this following (1.9). Let $r$ be an integer and consider the grading of $A_{\mathbf{Z}}$ defined by the rule $\mathrm{weight}(a_i) = i - s + r$ if $i \geqslant s$ and $\mathrm{weight}(a_i) = 0$ otherwise, then $\mathrm{redDisc}(h)$ is homogeneous of degree $(d-s+2r)(d-s-1)$. In particular, if $r = s$, i.e., weight$(a_i) = i$ if $i \geqslant s$ and 0 otherwise, we get that redDisc$(h)$ is homogeneous of degree

(2.11.1) $$(d+s)(d-s-1) = d(d-1) - s(s+1).$$

Similarly, if the grading of $A_{\mathbf{Z}}$ is defined by the rule weight$(a_i) = d - i + r$ if $i \geqslant s$ and 0 otherwise, then redDisc$(h)$ is homogeneous of the same degree $(d-s+2r)(d-s-1)$.

The following result is similar to Theorem (1.15). It is the key to the generalized Salmon formula for the discriminant.

**(2.12) Theorem** ([3, Theorem 4.25]). *Suppose that the ring $A_{\mathbf{Z}}$ is graded by means of the Zariski grading (1.14), i.e., weight$(u_{\boldsymbol{\alpha}}) = \max(\alpha_0 - d + s, 0)$. Then $\mathrm{Disc}(f)$ is of valuation $s(s-1)^n$ and its homogeneous part $H$ in this degree satisfies the following equality in $A_{\mathbf{Z}}$:*

(2.12.1) $$H \cdot \mathrm{Disc}(f_d) = \mathrm{Disc}(g) \cdot \mathrm{Disc}(f_s) \cdot \mathrm{redRes}(\partial_1 h, \ldots, \partial_n h, h).$$

Note that the three elements $\mathrm{Disc}(f_d)$, $\mathrm{Disc}(f_s)$, and $\mathrm{redRes}(\partial_1 h, \ldots, \partial_n h, h)$ are of degree 0 with respect to the Zariski grading, whereas $\mathrm{Disc}(g)$ is homogeneous of degree $s(s-1)^n$ by (2.5.2) and (1.8); recall that $f_d$ and $f_s$ (resp. $g$) are generic homogeneous polynomial of degrees $d$ and $s$ (resp. $s$) in the variables $(x_1, \ldots, x_n)$ (resp. $(x_0, \ldots, x_n)$).

**(2.13) Corollary.** *Using the Zariski grading of $A_{\mathbf{Z}}$ as in Theorem (2.12), the discriminant $\mathrm{Disc}(f)$ is of valuation $s(s-1)^n$ and can be written as*

$$\mathrm{Disc}(f) = \mathrm{Disc}(g) \cdot \mathrm{Disc}(f_s)^2 \cdot \mathrm{redDisc}(h) + \bigl(\text{terms of Zariski weight} > s(s-1)^n\bigr).$$

**(2.14)** This corollary provides an interesting connection between classical and reduced discriminants. As a first illustration of its interest, we give the following weight property of the reduced discriminant, which generalizes to arbitrary $n$ the computation of (2.11.1). Using the grading of $A_{\mathbf{Z}}$ defined by the rule weight$(u_\alpha) = d - \alpha_0$ as in (2.6), i.e., we give weight $j$ to all the coefficients of $f_j$, then the reduced discriminant redDisc$(h)$ is homogeneous of degree

$$n\bigl[d(d-1)^n - s(s+1)(s-1)^{n-1}\bigr].$$

Indeed, we know by (2.6) that $\mathrm{Disc}(f)$ is homogeneous of degree $nd(d-1)^n$ and that $\mathrm{Disc}(g)$ is homogeneous of degree $ns(s-1)^n$ by the same computation. In addition, all the coefficients of $f_s$ have weight $s$ in this grading and hence, applying Definition (2.2) we get that $\mathrm{Disc}(f_s)$ is homogeneous of degree $ns(s-1)^{n-1}$ (note that $f$ is a homogeneous polynomial in $n$ variables only). >From here the conclusion follows by a straightforward computation.

**(2.15) Comment.** In the above paragraphs we have introduced the reduced discriminant and provided some first properties that are sufficient for our purposes. However, a more detailed and complete study of this new eliminant polynomial, including for instance its irreducibility and the geometric meaning of its vanishing, in particular its connection to Milnor number, seems worthwhile and is left for future work.

## 2.3 – Application to the Salmon formula

We shall now see that the Salmon formula () is a particular case of the decomposition formula given in Corollary (2.13). We will then be able to generalize it to the case of a hypersurface in arbitrary dimension.

**(2.16)** Let $f(x, y, z)$ be the generic homogeneous polynomial of degree $d \geqslant 3$ and set

$$f = z^d f_0 + z^{d-1} f_1 + z^{d-2} f_2 + \cdots + z f_{d-1} + f_d,$$

$$f_0 = U, \qquad f_1 = Sx + Ty, \qquad f_2 = \tfrac{1}{2}(Ax^2 + 2Bxy + Cy^2),$$

(we introduce a rational number in the definition of $f_2$ only to follow Salmon's notation).

We consider the truncation of $f$ at order $d - 2$ with respect to $z$ and the corresponding Zariski grading. So, $U$ has degree 2, $S, T$ have degree 1, and the coefficients of the other $f_k$'s, $k \geqslant 2$, have degree 0 (this includes the coefficients $A, B, C$). We let

$$h = z^{d-2} f_2 + \cdots + z f_{d-1} + f_d \quad \text{and} \quad g = z^2 f_0 + z f_1 + f_2.$$

Then, Corollary (2.13) implies that

(2.16.1) $\qquad \mathrm{Disc}(f) = \mathrm{Disc}(g) \cdot \mathrm{Disc}(f_2)^2 \cdot \mathrm{redDisc}(h) + (\text{terms of Zariski weight} \geqslant 3).$

One has

$$\mathrm{Disc}(f_2) = \tfrac{1}{4} \begin{vmatrix} A & B \\ B & C \end{vmatrix} = \tfrac{1}{4}(AC - B^2)$$

and

$$\mathrm{Disc}(g) = \begin{vmatrix} A/2 & B/2 & S \\ B/2 & C/2 & T \\ S & T & U \end{vmatrix} = \frac{1}{4}\left(ACU - 2AT^2 - B^2 U + 4BST - 2CS^2\right)$$

so that (2.16.1) reads

$$\mathrm{Disc}(f) = \frac{1}{2^4}\left(ACU - 2AT^2 - B^2 U + 4BST - 2CS^2\right) \cdot (AC - B^2)^2 \cdot \mathrm{redDisc}(h)$$

$$\mod \left((S, T)^3 + U(S, T) + (U^2)\right).$$

Then, the specialization $U = S = 0$ yields the Salmon formula (🖋) with $\varphi = \mathrm{redDisc}(h)$; the latter is a homogeneous polynomial of degree $3(d-1)^2 - 7$ (for the standard grading). We mention that computations in the cases $d = 2$ and $d = 3$ with a computer algebra system have shown that $\mathrm{redDisc}(h)$ is an irreducible polynomial in these two cases.

**(2.17)** The Salmon formula for the discriminant of a plane curve can be generalized to the case of a hypersurface in a projective space of arbitrary dimension as follows. In a suitable system of homogeneous coordinates, any hypersurface $V(f) \subset \mathbf{P}^n$ has an equation of the form

(2.17.1) $\qquad f(x_0, x_1, \ldots, x_n) = T x_0^{d-1} x_n + \sum_{k=2}^{d} x_0^{d-k} f_k(x_1, \ldots, x_n)$

where for $k = 2, \ldots, d$ the polynomial $f_k$ is homogeneous of degree $k$ in the variables $x_1, \ldots, x_n$, respectively. We are here merely imposing that the hypersurface $V(f)$ goes through the point $(1 : 0 : \cdots : 0)$ and that its tangent hyperplane at this point is given by $x_n = 0$. Applying Corollary (2.13) as above, and setting $h = \sum_{k=2}^{d} x_0^{d-k} f_k$, we deduce that

$$\mathrm{Disc}(f) = \mathrm{Disc}(T x_0 x_n + f_2) \, \mathrm{Disc}(f_2)^2 \, \mathrm{redDisc}(h) \mod T^3.$$

Now, let $\bar{f}_2(x_1, \ldots, x_{n-1})$ be the homogeneous polynomial of degree 2 in the variables $x_1, \ldots, x_{n-1}$ defined as $\bar{f}_2 = f_2(x_1, \ldots, x_{n-1}, 0)$ (beware the difference in notation with (2.5)). Then, we get

$$
\begin{aligned}
2^{a(n+1,2)} \operatorname{Disc}(Tx_0 x_n + f_2) &= \operatorname{Res}\big(\partial_0(Tx_0 x_n + f_2), \partial_1(Tx_0 x_n + f_2), \ldots, \partial_n(Tx_0 x_n + f_2)\big) \\
&= T \cdot \operatorname{Res}\big(\partial_1 f_{2|x_n=0}, \ldots, \partial_{n-1} f_{2|x_n=0}, Tx_0 + \partial_n f_{2|x_n=0}\big) \\
&= T^2 \cdot \operatorname{Res}(\partial_1 \bar{f}_2, \ldots, \partial_{n-1} \bar{f}_2) = 2^{a(n-1,2)} T^2 \operatorname{Disc}(\bar{f}_2).
\end{aligned}
$$

It follows that $\operatorname{Disc}(Tx_0 x_n + f_2) = T^2 \operatorname{Disc}(\bar{f}_2)$, and eventually we obtain the following generalized Salmon formula:

(2.17.2) $$\operatorname{Disc}(f) = T^2 \operatorname{Disc}(\bar{f}_2) \operatorname{Disc}(f_2)^2 \operatorname{redDisc}(h) + T^3 \psi.$$

# 3 – Computation of the node-couple degree by elimination

In [16, § 605–607], Salmon sets up the following strategy to compute the number of 2-nodal curves in a general net of hyperplane sections of a smooth (hyper)surface $S \subset \mathbf{P}^3$. For $p' \in S$ and $p'' \in \mathbf{T}_{p'} S - p'$, consider the pencil $\mathfrak{l}_{p',p''} \subset \check{\mathbf{P}}^3$ of (hyper)planes containing $p'$ and $p''$ (i.e., $\mathfrak{l}_{p',p''} = \langle p', p'' \rangle^\perp$). It cuts out $d^\vee = \deg(S^\vee)$ points on the dual surface $S^\vee$, counted with multiplicities, corresponding to planes tangent to $S$. Among these, $\mathbf{T}_{p'} S$ counts doubly if it is a plain tangent plane, and triply if it is plainly bitangent.

Indeed, the line $\mathfrak{l}_{p',p''}$ is contained in the plane $(p')^\perp$, hence tangent to $S^\vee$ at the point $(\mathbf{T}_{p'} S)^\perp \in S^\vee$. If the plane $\mathbf{T}_{p'} S$ is plainly bitangent to $S$, then $(\mathbf{T}_{p'} S)^\perp$ is a general point on the ordinary double curve of $S^\vee$, and the line $\mathfrak{l}_{p',p''}$ is tangent to one of the two transverse sheets of $S^\vee$ at $(\mathbf{T}_{p'} S)^\perp$. The idea is then firstly to determine the conditions on $p'$ for $\mathbf{T}_{p'} S$ to count with multiplicity greater than 2 in $\mathfrak{l}_{p',p''} \cap S^\vee$, and secondly to sort out the various corresponding geometric situations. A key element to carry this out is the famous formula (🕭); another one is the elimination procedure (3.12).

We work out Salmon's procedure in subsection 3.2, and in subsection 3.3 we show how it carries over for hypersurfaces in a projective space of arbitrary dimension. >From now on, we work over an algebraically closed field $\mathbf{k}$ of characteristic 0.

## 3.1 – Polarity

We give here a brief recap on polarity, so that the reader unfamiliar with this may conveniently consult the relevant material. We refer the reader to [8, Chapter 1] and [1, § 5.4 and 5.6] for more details.

**(3.1) The polar pairing.** Let $f$ be a complex homogeneous degree $d$ polynomial in $n+1$ variables $(x_0, \ldots, x_n)$. For $\mathbf{a} = (a_0, \ldots, a_n) \in \mathbf{C}^{n+1}$, we define

(3.1.1) $$\mathrm{D}_a f = a_0 \partial_0 f + \cdots a_n \partial_n f,$$

and for $\mathbf{a}_1, \ldots, \mathbf{a}_k \in \mathbf{C}^{n+1}$, we let

(3.1.2) $$\mathrm{D}_{\mathbf{a}_1 \cdots \mathbf{a}_k} f = (\mathrm{D}_{\mathbf{a}_1} \circ \cdots \circ \mathrm{D}_{\mathbf{a}_k})(f).$$

Extending (3.1.2) by $\mathbf{C}$-linearity, one obtains a perfect bilinear pairing

$$\mathrm{D} : \operatorname{Sym}^k(\mathbf{C}^{n+1}) \otimes \operatorname{Sym}^d(\check{\mathbf{C}}^{n+1}) \longrightarrow \operatorname{Sym}^{d-k}(\check{\mathbf{C}}^{n+1})$$
$$\xi \otimes f \longmapsto \mathrm{D}_\xi f$$

extending the natural pairing between $\mathbf{C}^{n+1}$ and $\check{\mathbf{C}}^{n+1}$.

**(3.2) Polar hypersurfaces.** Consider the hypersurface $X = V(f) \subset \mathbf{P}^n$. For $\hat{a} = (a_0, \ldots, a_n) \in \mathbf{C}^{n+1}$, and $k \in \mathbf{N}$,
$$\mathrm{D}_{\hat{a}^k} f = \bigl(a_0 \partial_0 + \cdots + a_n \partial_n\bigr)^k f,$$
may be viewed as a bihomogeneous polynomial of bidegree $(k, d-k)$ in the variables $(a_0, \ldots, a_n)$ and $(x_0, \ldots, x_n)$. The hypersurface $V(\mathrm{D}_{\hat{a}^k} f)$ depends only on $X$ and the point $a = (a_0 : \ldots : a_n) \in \mathbf{P}^n$; we call it the *k-th polar of $X$ with respect to $a$*, and denote it by $\mathrm{D}_{a^k} X$. We will often abuse notation and consider $\hat{a}$ and $a$ without distinction.

We shall also use the following useful notation: $\mathrm{D}^k X(a)$, referred to as the *polar k-ic of $X$ at $a$*, is the hypersurface defined by the degree $k$ polynomial in the variables $(x_0, \ldots, x_n)$,
$$\mathrm{D}^k f(\hat{a}) = \bigl(x_0 \partial_0 + \cdots x_n \partial_n\bigr)^k f(a_0, \ldots, a_n).$$

**(3.3) Proposition** (polar symmetry). *Consider $X \subset \mathbf{P}^n$ a degree $d$ hypersurface. Let $a, b \in \mathbf{P}^n$, and $k \in [\![1, d-1]\!]$. One has the equivalence:*
$$a \in \mathrm{D}_{b^k} X \iff b \in \mathrm{D}_{a^{d-k}} X.$$

This says that $\mathrm{D}^k X(a) = \mathrm{D}_{a^{d-k}} X$. In case $a$ is a singular point of $\mathrm{D}_{b^k} X$, we have the following. Let $k, l$ be positive integers such that $k + l < d$. If $a$ is a point of multiplicity $\geqslant d - k - l + 1$ of $\mathrm{D}_{b^l} X$, then $b$ is a point of multiplicity $\geqslant d - k - l + 1$ of $\mathrm{D}_{a^k} X$. The proof is mere polynomial calculus.

**(3.4)** If $b \in X = V(f)$ is a smooth point, then $\mathrm{D}^1 f(b)$ is "the" linear homogeneous polynomial defining the tangent hyperplane to $X$ at $b$. Therefore, for all $a \in \mathbf{P}^n$,
$$X \cap \mathrm{D}_a X = \{x \in X : \mathbf{T}_x X \ni a\}.$$
This generalizes to the following fundamental property.

If $X \subset \mathbf{P}^n$ is a hypersurface and $\ell \subset \mathbf{P}^n$ a line, for all $p \in \ell$ we let $i(X, \ell)_p$ be the multiplicity with which $p$ appears in $X \cap \ell$.

**(3.5) Theorem.** *Let $X$ be a degree $d$ hypersurface, $a \in X$, and $b \in \mathbf{P}^n$. For all integer $s \geqslant 0$, one has*
$$i\bigl(X, \langle a, b \rangle\bigr)_a \geqslant s+1 \iff \forall k \leqslant s, \quad b \in \mathrm{D}_{a^{d-k}} X$$
$$\iff \forall k \leqslant s, \quad a \in \mathrm{D}_{b^k} X.$$

It turns out that for $a \in X$, all polars of $X$ with respect to $a$ (equivalently, all $\mathrm{D}^k X(a)$) are tangent at $a$. Actually, $X$ and its polar $k$-ic at $a$, $\mathrm{D}^k X(a)$ have the same polar $s$-ics at $a$ for all $s \leqslant k$, as the following identities show:
$$\mathrm{D}^s \bigl(\mathrm{D}^k X(a)\bigr)(a) = \mathrm{D}_{a^{k-s}}\bigl(\mathrm{D}^k X(a)\bigr) = \mathrm{D}_{a^{k-s}}\bigl(\mathrm{D}_{a^{d-k}} X\bigr) = \mathrm{D}_{a^{d-s}} X = \mathrm{D}^s X(a).$$
This has the following remarkable consequence.

**(3.6) Corollary.** *Let $X$ be a hypersurface in $\mathbf{P}^n$, and $a$ a point of $X$. The $n$ hypersurfaces $\mathrm{D}^1 X(a), \ldots, \mathrm{D}^n X(a)$ intersect with multiplicity at least $n!$ at $a$.*

*Proof.* We consider only the case in which all intersections are complete, otherwise the result is trivial. Then the intersection $\mathrm{D}^1 X(a), \ldots, \mathrm{D}^{n-1} X(a)$ consists of $(n-1)!$ lines each intersecting $X$ with multiplicity at least $n$ at $a$, as follows from Theorem (3.5). Since the polar $k$-ics at $a$ of $D^n X(a)$ are the same as those of $X$ as indicated above, each of these lines also intersect $D^n X(a)$ with multiplicity at least $n$ at $a$, and the result follows. $\square$

The polar hyperplane $D^1 X(a)$ is well-defined only if $a$ is a smooth point of $X$ (otherwise its equation is 0, and the more reasonable thing to do is to set $D^1 X(a) = \mathbf{P}^n$). When $a$ is singular, the following holds.

**(3.7) Theorem.** *Let $X$ be a degree $d$ hypersurface, and $a \in X$ a point of multiplicity $m$. We consider an integer $r \leqslant d - m$.*
*(3.7.1) The polar hypersurface $D_{a^r} X$ has multiplicity $m$ at $a$, and it has the same tangent cone at $a$ as $X$: $TC_a(D_{a^r} X) = TC_a(X)$.*
*(3.7.2) Let $b \in \mathbf{P}^n - \{a\}$. The polar hypersurface $D_{b^r} X$ has multiplicity $\geqslant m - r$ at $a$; this multiplicity is exactly $m - r$ for general $b \in \mathbf{P}^n$, and in this case the tangent cone of $D_{b^r} X$ at $a$ equals the $r$-th polar with respect to $b$ of the tangent cone $TC_a(X)$:*

$$TC_a(D_{b^r} X) = D_{b^r}(TC_a(X)).$$

In particular, (3.7.2) tells us that $D_{b^r} X$ contains $a$ if $r \leqslant m - 1$, and is singular at $a$ if $r \leqslant m - 2$.

We end this section by recalling the following definition.

**(3.8) Hessian of a hypersurface.** Let $X$ be a hypersurface in $\mathbf{P}^n$, defined by the homogeneous polynomial $f \in \mathbf{C}[x_0, \ldots, x_n]$. The *Hessian hypersurface* of $X$ is the hypersurface defined by the homogeneous polynomial

$$\text{Hess}(f) = \det(\partial_i \partial_j f)_{0 \leqslant i,j \leqslant n}.$$

## 3.2 – The computation

**(3.9) Setup.** We recall that $\mathbf{k}$ is assumed to be an algebraically closed field of characteristic 0. We consider $S \subset \mathbf{P}^3$ a smooth surface of degree $d$, defined by a homogeneous polynomial $f(x, y, z, w)$. Let $\hat{p}' = (x', y', z', w')$, $\hat{p}'' = (x'', y'', z'', w'')$, $\hat{p} = (x, y, z, w) \in \mathbf{k}^4$, and call $p', p'', p$ the corresponding points in $\mathbf{P}^3$ (beware the unusual distribution of the prime decorations). The choice of $\hat{p}', \hat{p}'', \hat{p}$ defines a system of homogeneous coordinates $(\alpha : \beta : \gamma)$ on the plane generated by $p', p''$ and $p$. One has
(3.9.1)

$$\begin{aligned} f(\alpha\hat{p}' + \beta\hat{p}'' + \gamma\hat{p}) &= \alpha^d \cdot f(\hat{p}') + \alpha^{d-1} \cdot D_{\beta\hat{p}''+\gamma\hat{p}} f(\hat{p}') + \frac{1}{2}\alpha^{d-2} \cdot D_{(\beta\hat{p}''+\gamma\hat{p})^2} f(\hat{p}') \mod (\beta,\gamma)^3 \\ &= \alpha^d \cdot f(\hat{p}') + \alpha^{d-1}\big(\beta D_{\hat{p}''} f(\hat{p}') + \gamma D_{\hat{p}} f(\hat{p}')\big) \\ &\quad + \frac{1}{2}\alpha^{d-2}\big(\beta^2 D_{\hat{p}''^2} f(\hat{p}') + 2\beta\gamma D_{\hat{p}''\hat{p}} f(\hat{p}') + \gamma^2 D_{\hat{p}^2} f(\hat{p}')\big) \mod (\beta,\gamma)^3. \end{aligned}$$

Considered as a homogeneous polynomial in the variables $(\alpha, \beta, \gamma)$, this is the equation of the hyperplane section of $S$ by $\langle p', p'', p \rangle$.

Choose $p' \in S$ and $p'' \in \mathbf{T}_{p'} S$, and think of them as fixed for a moment. Then $f(\hat{p}') = D_{\hat{p}''} f(\hat{p}') = 0$, so that (3.9.1) reduces to

(3.9.2) $\qquad f(\alpha\hat{p}' + \beta\hat{p}'' + \gamma\hat{p}) = T\alpha^{d-1}\gamma + \frac{1}{2}\alpha^{d-2}(A\beta^2 + 2B\beta\gamma + C\gamma^2) \mod (\beta,\gamma)^3,$

where $T = D_{\hat{p}}f(\hat{p}')$, $A = D_{\hat{p}''^2}f(\hat{p}')$, $B = D_{\hat{p}''\hat{p}}f(\hat{p}')$, $C = D_{\hat{p}^2}f(\hat{p}')$. Note in particular that, as homogeneous polynomials in the variables $\hat{p} = (x, y, z, w)$, $T$ and $C$ are the equations of the tangent plane $\mathbf{T}_{p'}S$ and the polar quadric $D^2S(p')$ of $S$ at $p'$ respectively.

Let us consider the discriminant of the plane curve $S \cap \langle p', p'', p \rangle$, i.e., the discriminant of $f(\alpha \hat{p}' + \beta \hat{p}'' + \gamma \hat{p})$ as a homogeneous polynomial in $(\alpha, \beta, \gamma)$; by (▲), it writes

$$\text{(3.9.3)} \qquad \text{Disc}\bigl(f(\alpha \hat{p}' + \beta \hat{p}'' + \gamma \hat{p})\bigr) = T^2\bigl(A(B^2 - AC)^2 \varphi + T\psi\bigr).$$

It vanishes if and only if the plane $\langle p', p'', p \rangle$ is tangent to $S$. It is a trihomogeneous polynomial in the variables $\hat{p}'$, $\hat{p}''$ and $\hat{p}$; considering $p'$ and $p''$ as fixed, we find a homogeneous polynomial in $\hat{p} = (x, y, z, w)$ vanishing along the $d^\vee$ planes (counted with multiplicities) tangent to $S$ in the pencil $\mathfrak{l}_{p',p''}$, with the notation of the introduction to the present section.

The fact that $T^2$ factors out of this discriminant gives an algebraic proof of the fact that the tangent plane $\mathbf{T}_{p'}S$ always appears with multiplicity $\geq 2$ in the scheme $\mathfrak{l}_{p',p''} \cap S^\vee$. We shall now derive the conditions under which it appears with multiplicity $> 2$, equivalently $T$ divides $A(B^2 - AC)^2 \varphi$.

Let us first give the key technical tools in elimination theory needed to carry this out. We begin with the following characterization of the non-emptiness of the intersection of two hyperplanes and a line in $\mathbf{P}^3$.

**(3.10) Lemma.** *Let $\hat{a}', \hat{a}'' \in \mathbf{k}^4$, and consider the line $L$ defined by the two linear homogeneous equations ${}^\mathsf{T}\hat{a}'(\hat{p}) = {}^\mathsf{T}\hat{a}''(\hat{p}) = 0$. The intersection of the two lines $L$ and $\langle p', p'' \rangle$ is non-empty if and only if the irreducible polynomial*

$$D(a', a'', p', p'') := \begin{vmatrix} {}^\mathsf{T}\hat{a}'(\hat{p}') & {}^\mathsf{T}\hat{a}''(\hat{p}') \\ {}^\mathsf{T}\hat{a}'(\hat{p}'') & {}^\mathsf{T}\hat{a}''(\hat{p}'') \end{vmatrix} = {}^\mathsf{T}\hat{a}'(\hat{p}') \cdot {}^\mathsf{T}\hat{a}''(\hat{p}'') - {}^\mathsf{T}\hat{a}'(\hat{p}'') \cdot {}^\mathsf{T}\hat{a}''(\hat{p}')$$

*vanishes.*

*Proof.* The line $\langle p', p'' \rangle$ is the image of the map $\mathbf{P}^1 \to \mathbf{P}^3$ that sends the point $(u : v) \in \mathbf{P}^1$ to the point $u\hat{p}' + v\hat{p}'' \in \mathbf{P}^3$. Therefore, the lines $L$ and $\langle p', p'' \rangle$ intersect if and only if the two polynomials ${}^\mathsf{T}\hat{a}'(u\hat{p}' + v\hat{p}'')$ and ${}^\mathsf{T}\hat{a}''(u\hat{p}' + v\hat{p}'')$ share a common root in $\mathbf{P}^1$. As these polynomials are linear forms in $u, v$, the polynomial $D(a', a'', p', p'')$ is simply the determinant of the corresponding linear system.

It remains to show that $D$ is an irreducible polynomial. Observe that it is a linear form in each set of variables $a'$, $a''$, $p'$ and $p''$. Set $\hat{a}' = (a'_1, \ldots, a'_4) \in \mathbf{k}^4$ and take similar notations for $a'', p'$ and $p''$; we can write

$$D(a', a'', p', p'') = \sum_{i=1}^4 a'_i \begin{vmatrix} p'_i & {}^\mathsf{T}\hat{a}''(\hat{p}') \\ p''_i & {}^\mathsf{T}\hat{a}''(\hat{p}'') \end{vmatrix} = \sum_{i=1}^4 a'_i D_i(a'', p', p'').$$

If $D$ factors as a product $P.Q$, then $P$ or $Q$, say $Q$, must be independent of $a'$ and $P$ must be a linear form in $a'$:

$$P = a'_1 P_1(a'', p', p'') + \cdots + a'_4 P_4(a'', p', p'').$$

It follows that we must have $QP_i = D_i$ for all $i = 1, \ldots, 4$. Now, we claim that $D_1, D_2, D_3$ and $D_4$ are irreducible and coprime polynomials, which implies that $Q$ must be equal to $\pm 1$ and proves the irreducibility of $D$.

To prove the claim, we proceed similarly since $D_1, \ldots, D_4$ are linear forms in $a'', p'$ and $p''$. For all $i = 1, \ldots, 4$ we have

$$(3.10.1) \qquad D_i(a'', p', p'') = \sum_{j=1}^{4} a_j'' \begin{vmatrix} p_i' & p_j' \\ p_i'' & p_j'' \end{vmatrix}.$$

We deduce that $D_i$ must be irreducible for otherwise it would have a factor that is independent of $a''$, but then this factor must be a common factor of the minors appearing in (3.10.1), which is impossible because these minors are known to be irreducible and coprime polynomials. In addition, the coprimeness of $D_1, \ldots, D_4$ follows in a straightforward manner. □

**(3.11) Corollary.** *Let $\hat{a}', \hat{a}'' \in \mathbf{k}^4$, and consider in addition two linear forms $H_1(\hat{p})$, $H_2(\hat{p})$, with coefficients depending on $\hat{p}', \hat{p}''$, such that*

$$\forall \hat{p}', \hat{p}'', \quad H_1(\hat{p}') = H_1(\hat{p}'') = H_2(\hat{p}') = H_2(\hat{p}'') = 0.$$

*Then the resultant of $H_1, H_2, {}^\mathsf{T}\hat{a}', {}^\mathsf{T}\hat{a}''$, with respect to the homogeneous variable $\hat{p}$, is divisible by the polynomial $D(a', a'', p', p'')$ defined in Lemma (3.10).*

*Proof.* The equations $H_1$ and $H_2$ define two hyperplanes that both contain the two points $p'$ and $p''$, hence the whole line $\langle p', p'' \rangle$. Therefore, if there exists an intersection point between the lines $\langle p', p'' \rangle$ and $V({}^\mathsf{T}\hat{a}', {}^\mathsf{T}\hat{a}'')$, then the hyperplanes $H_1, H_2, {}^\mathsf{T}\hat{a}'$ and ${}^\mathsf{T}\hat{a}''$ intersect at this point and hence their resultant vanishes. Using Lemma (3.10), the claimed property follows. □

Now, let $W_{\hat{p}', \hat{p}''}(\hat{p})$ be a trihomogeneous polynomial of tridegree $(\lambda, \mu, \mu)$ in the sets of variables $(\hat{p}', \hat{p}'', \hat{p})$, for some non-negative integers $\lambda$ and $\mu$, and assume that for all

$$(p', p'') \in \tilde{\mathcal{C}}_S := \{(x', x'') : x' \in S \text{ and } x'' \in \mathbf{T}_{x'}S\} \subset \mathbf{P}^3 \times \mathbf{P}^3$$

the hypersurface $V(W_{\hat{p}', \hat{p}''}) \subset \mathbf{P}^3$ consists of $\mu$ planes, counted with multiplicities, all containing the line $\langle p', p'' \rangle$.

**(3.12) Proposition.** *For any $p' \in S$, there exists a $p'' \in \mathbf{T}_{p'}S - \{p'\}$ such that the tangent plane $\mathbf{T}_{p'}S$ is a component of $V(W_{\hat{p}', \hat{p}''})$ if and only if the same holds for all $p'' \in \mathbf{T}_{p'}S - \{p'\}$. The set of those $p' \in S$ for which this happens is cut out on $S$ by a hypersurface of degree $\lambda + \mu(\deg(S) - 2)$.*

*Proof.* The idea is to express in terms of resultants the trivial fact that $\mathbf{T}_{p'}S$ is a component of $V(W_{\hat{p}', \hat{p}''})$ if and only if for all line $L \subset \mathbf{P}^3$ the intersection $L \cap \mathbf{T}_{p'}S \cap V(W_{\hat{p}', \hat{p}''})$ is non-empty.

Consider $T = \mathrm{D}_{\hat{p}} f(\hat{p}')$ as a trihomogeneous polynomial of tridegree $(d-1, 0, 1)$ in $(\hat{p}', \hat{p}'', \hat{p})$. Let $\hat{a}', \hat{a}'' \in \mathbf{k}^4$, and consider the line $L$ defined by the two linear homogeneous equations

$${}^\mathsf{T}\hat{a}'(\hat{p}) = {}^\mathsf{T}\hat{a}''(\hat{p}) = 0.$$

It intersects $\mathbf{T}_{p'}S \cap V(W_{\hat{p}', \hat{p}''})$ if and only if the resultant $R = \mathrm{Res}(T, W, {}^\mathsf{T}\hat{a}', {}^\mathsf{T}\hat{a}'')$, with respect to the variables $\hat{p}$, vanishes. It follows from the homogeneity properties of the resultant (see Theorem (1.6)) that $R$ is a multi-homogeneous polynomial in the variables $(\hat{a}', \hat{a}'', \hat{p}', \hat{p}'')$, of multi-degree $(\mu, \mu, \mu(d-1) + \lambda, \mu)$.

For $(p', p'') \in \tilde{\mathcal{C}}_S$, the polynomial $W$ splits as the product of $\mu$ linear forms $W_i$, $i = 1, \ldots, \mu$ that all contain the line $\langle p', p'' \rangle$. Therefore, by the divisibility property of the resultant, $R$ as well splits as the product of the resultant $\mathrm{Res}(T, W_i, {}^\mathsf{T}\hat{a}', {}^\mathsf{T}\hat{a}'')$. Therefore, applying Corollary (3.11) we deduce that there exists a multi-homogeneous polynomial $R'$ such that after restriction to $\tilde{\mathcal{C}}_S \times \mathbf{P}^3 \times \mathbf{P}^3$ one has $R = D(a', a'', p', p'')^\mu R'$. Computing degrees one finds that $R'$ has multi-degree $(0, 0, \mu(d-2) + \lambda, 0)$ in the variables $(\hat{a}', \hat{a}'', \hat{p}', \hat{p}'')$, which proves the statement. □

We now discuss the various possibilities for $T$ to divide $A(B^2 - AC)^2 \varphi$, in the setup of (3.9).

**(3.13) Vanishing of $A$.** The polynomial $A$ has degree 0 in $p$, therefore it is divisible by $T$ if and only if it is identically zero. By definition $A = \mathrm{D}_{\hat{p}''^2} f(\hat{p}')$, so its vanishing is equivalent to the point $p''$ being on the polar quadric $\mathrm{D}^2 S(p')$. Since $p'' \in \mathbf{T}_{p'} S$, this in turn is equivalent to the line $\langle p', p'' \rangle$ being one of the two inflectional tangents of $S$ at $p'$. We thus find an algebraic proof of the fact that if $i(\langle p', p'' \rangle, S)_{p'} \geqslant 3$, then $i(\mathfrak{l}_{p', p''}, S^\vee)_{(\mathbf{T}_{p'} S)^\perp} \geqslant 3$, which is a manifestation of biduality.

Geometrically, if a line $L \subset \mathbf{P}^3$ is an inflexional tangent to $S$ at $p'$ (i.e., $L$ is the tangent line to one of the local branches of $\mathbf{T}_p S \cap S$ at $p'$), then its orthogonal $L^\perp \subset \check{\mathbf{P}}^3$ is an inflexional tangent to $S^\vee$ at the point $(\mathbf{T}_{p'} S)^\perp$.

As a side remark note that when $A = 0$, all the curves cut out on $S$ by a member of the pencil $\mathfrak{l}_{p', p''}$ have an inflexion point at $p'$.

**(3.14)** The divisibilities of $B^2 - AC$ and $\varphi$ by $T$ are analyzed using Proposition (3.12). To see that the latter result indeed applies, we note that for $(p', p'') \in \tilde{\mathcal{C}}_S$ the discriminant (3.9.3), $T^2 A(B^2 - AC)^2 \varphi + T^3 \psi$, defines as a homogeneous polynomial in the variable $p$ a hypersurface consisting of $d^\vee$ planes, counted with multiplicities, all containing the line $\langle p', p'' \rangle$. This implies that so does its homogeneous part of lowest degree with respect to the Zariski grading, namely $T^2 A(B^2 - AC)^2 \varphi$ (see subsection 2.3), hence also all of its factors. This is obvious for $T$; the polynomial $A$ on the other hand is independent on the variable $p$, hence defines either the whole space, or the empty set.

**(3.15) Divisibility of $B^2 - AC$ by $T$.** The polynomial $B^2 - AC$ is tri-homogeneous of tri-degree $(2(d-2), 2, 2)$ in the variables $(p', p'', p)$. It follows from Proposition (3.12) together with (3.14) that there exists a homogeneous polynomial $H$ of degree $4(d-2)$ in the variable $p'$, with constant coefficients, such that for fixed $p' \in S$ and $p'' \in \mathbf{T}_{p'} S$, $T$ divides $B^2 - AC$ as polynomials in $p$ if and only if $H(p') = 0$.

We recover the Hessian determinant. Indeed, it follows from (3.9.2) that the tangent cone of the section of $S$ by its tangent hyperplane at $p'$, $S \cap \mathbf{T}_{p'} S$, is given by the symmetric matrix

$$\begin{pmatrix} A & B \\ B & C \end{pmatrix}.$$

Thus, $B^2 - AC$ is zero modulo $T$ if and only if the curve $S \cap \mathbf{T}_{p'} S$ has a degenerate tangent cone at $p'$, equivalently $p'$ is a parabolic point of $S$.

Geometrically this is explained as follows. If $p'$ is a parabolic point of $S$, then the point $(\mathbf{T}_{p'} S)^\perp \in \check{\mathbf{P}}^3$ sits on the cuspidal double curve of $S^\vee$, and for general $p'' \in \mathbf{T}_{p'} S$ the line $\mathfrak{l}_{p', p''} \subset \check{\mathbf{P}}^3$ is contained in the tangent cone of $S^\vee$ at $(\mathbf{T}_{p'} S)^\perp \in \check{\mathbf{P}}^3$, hence $i(\mathfrak{l}_{p', p''}, S^\vee)_{(\mathbf{T}_{p'} S)^\perp} \geqslant 3$.

**(3.16) Divisibility of the reduced discriminant $\varphi$ by $T$.** The discriminant (3.9.3) is tri-homogeneous in the variables $(p', p'', p)$, of degree $d^\vee = d(d-1)^2$ with respect to all three variables. Indeed, if we fix two of the three points $(p', p'', p)$, $p_1$ and $p_2$ say, then we get a polynomial whose zero locus in $\mathbf{P}^3$ is the sum of all hyperplanes tangent to $S$ in the pencil defined by the line $\langle p_1, p_2 \rangle$, and there are $\deg(S^\vee) = d^\vee$ of them, counted with mutiplicities.

On the other hand, the fact that (3.9.3) is of tri-degree $\big(d(d-1)^2, d(d-1)^2, d(d-1)^2\big)$ is a straightforward application of (2.6); this proves the equality $d^\vee = d(d-1)^2$.

This implies that also the homogeneous piece of (3.9.3) of lowest degree with respect to the Zariski grading is tri-homogeneous of tri-degree $(d^\vee, d^\vee, d^\vee)$ for the standard grading. Computing degrees, one then finds that $\varphi$ is of tri-degree $(\lambda, \mu, \mu)$ in $(p', p'', p)$, with $\lambda = (d-2)(d^2-6)$ and $\mu = d^3 - 2d^2 + d - 6$. It thus follows as in (3.15) from Proposition (3.12) together with (3.14) that there exists a homogeneous polynomial $K$ of degree $(d-2)(d^3 - d^2 + d - 12)$ in the variable $p'$, with constant coefficients, such that for fixed $p' \in S$ and $p'' \in \mathbf{T}_{p'}S$, $T$ divides $\varphi$ as polynomials in $p$ if and only if $K(p') = 0$.

In conclusion, we have proven the following statement.

**(3.17) Theorem.** *Let $S$ be a smooth degree $d$ surface in $\mathbf{P}^3$. There is a hypersurface $V(K)$ of degree $(d-2)(d^3 - d^2 + d - 12)$, the intersection of which with $S$ is the locus of tangency points of planes bitangent to $S$.*

**(3.18) Corollary.** *The ordinary double curve of the dual surface $S^\vee$ has degree $\frac{1}{2}d(d-1)(d-2)(d^3 - d^2 + d - 12)$.*

*Proof.* Let $p'' \in \mathbf{P}^3$ be a general point. The locus of those points $p' \in S$ such that there exists a plane through $p''$ tangent to $S$ at $p'$ is the apparent boundary $\mathrm{D}_{p''}S \cap S$. Therefore, by Theorem (3.17), the locus of points $p' \in S$ such that $\mathbf{T}_{p'}S$ is bitangent and passes through $p''$ is $\mathrm{D}_{p''}S \cap S \cap V(K)$. Now for each bitangent plane there are two tangency points $p'$, so the number of bitangent planes passing through $p''$ is

$$\frac{1}{2} \cdot \deg(S) \cdot \deg(\mathrm{D}_{p''}S) \cdot \deg K.$$

□

## 3.3 – Generalization to hypersurfaces of arbitrary dimension

It turns out that Salmon's strategy actually works in arbitrary dimension, using the generalization (2.17.2) of formula (🦉). The arguments are exactly the same, so we are going to be sketchy.

**(3.19)** Let $X = V(f)$ be a smooth hypersurface of degree $d$ in $\mathbf{P}^n$, with $f$ in normal form as in (2.17.1), and consider $n$ points $p_1, \ldots, p_{n-1}, p$ (or rather $\hat{p}_1, \ldots, \hat{p}_{n-1}, \hat{p} \in \mathbf{k}^{n+1}$). The points $p_1, \ldots, p_{n-1}$ define a pencil of hyperplanes in $\mathbf{P}^n$, and for each $p$ the vectors $\hat{p}_1, \ldots, \hat{p}_{n-1}, \hat{p}$ define a system of homogeneous coordinates $(\alpha_1 : \ldots : \alpha_{n-1} : \alpha)$ on the member $H_p$ determined by $p$ in the pencil.

One may express in this system of coordinates the equation of the hyperplane section $H_p \cap X$, viz.

$$\tilde{f}_{p_1, \ldots, p_{n-1}, p}(\alpha_1, \ldots, \alpha_{n-1}, \alpha) := f(\alpha_1 \hat{p}_1 + \cdots + \alpha_{n-1}\hat{p}_{n-1} + \alpha \hat{p})$$
$$= T\alpha_1^{d-1}\alpha + \frac{1}{2}\alpha_1^{d-2} \cdot \mathrm{D}_{(\alpha_2\hat{p}_2 + \cdots + \alpha\hat{p})^2}f(\hat{p}_1) \mod (\alpha_2, \ldots, \alpha)^3,$$

where $T = \mathrm{D}_{\hat{p}}f(\hat{p}_1)$, and apply (2.17.2) to obtain a Taylor expansion of its discriminant, viz.

$$\mathrm{Disc}(H_p \cap X) = T^2 \cdot \mathrm{Disc}(\bar{F}_2) \cdot \mathrm{Disc}(F_2)^2 \cdot \varphi \mod T^3,$$

with $F_2(\alpha_2, \ldots, \alpha_n, \alpha) = \mathrm{D}_{(\alpha_2\hat{p}_2 + \cdots + \alpha\hat{p})^2}f(\hat{p}_1)$ and $\bar{F}_2(\alpha_2, \ldots, \alpha_n) = F_2(\alpha_2, \ldots, \alpha_n, 0)$. Our task is to analyze the divisibilities of $\mathrm{Disc}(F_2)$ and $\varphi$ by $T$.

**(3.20)** Lemma (3.10) and its Corollary (3.11) generalize as follows. Given $n-1$ linear homogeneous equations corresponding to $a_1, \ldots, a_{n-1} \in \mathbf{P}^n$, the line $V({}^\mathsf{T}a_1, \ldots, {}^\mathsf{T}a_{n-1})$ intersects the $(n-2)$-dimensional linear space $\langle p_1, \ldots, p_{n-1}\rangle$ if and only if

$$D(a_1, \ldots, a_{n-1}, p_1, \ldots, p_{n-1}) := \det\bigl({}^\mathsf{T}a_i(p_j)\bigr)_{1 \leqslant i,j \leqslant n-1}$$

vanishes. Then the natural adaptation of Corollary (3.11) holds. The proofs are mutatis mutandis the same as those of Lemma (3.10) and Corollary (3.11).

**(3.21)** Next, using (3.20), it is straightforward to adapt Proposition (3.12) and its proof. The upshot is the following: let

$$\tilde{\mathcal{C}}_X := \{(x_1, \ldots, x_{n-1}) \in \mathbf{P}^n \times \mathbf{P}^n : x_1 \in S \text{ and } x_2, \ldots, x_{n-1} \in \mathbf{T}_{x'}S\},$$

and consider a multihomogeneous polynomial $W_{p_1, \ldots, p_{n-1}}(p)$ of multidegree $(\lambda, \mu, \ldots, \mu)$ in the variables $p_1, \ldots, p_{n-1}, p$, such that for all $(p_1, \ldots, p_{n-1}) \in \tilde{\mathcal{C}}_X$ the hypersurface $V(W_{p_1, \ldots, p_{n-1}})$ consists of $\mu$ hyperplanes, counted with multiplicities, that all contain $\langle p_1, \ldots, p_{n-1}\rangle$.

The locus of those $p_1 \in X$ such that there exist $p_2, \ldots, p_{n-1} \in \mathbf{T}_{p_1}X$ such that (equivalently, for all $p_2, \ldots, p_{n-1} \in \mathbf{T}_{p_1}X$) $\mathbf{T}_{p_1}X$ is a component of $V(W_{p_1, \ldots, p_{n-1}})$ is cut out on $X$ by a hypersurface of degree $\lambda + \mu(d-2)$.

**(3.22)** One may then argue as in (3.14) to analyze the divisibilities of $\mathrm{Disc}(F_2)$ and $\varphi$ by $T$. For $\mathrm{Disc}(F_2)$, whe shall find the Hessian determinant of $X$, since as before $F_2 \pmod{T}$ defines the tangent cone at $p_1$ of $\mathbf{T}_{p_1}X \cap X$, so $\mathrm{Disc}(F_2)$ vanishes modulo $T$ if and only if the latter tangent cone is singular, i.e., if and only if $p_1$ is a parabolic point of $X$.

And indeed $\mathrm{Disc}(F_2)$, being the determinant of the symmetric matrix

$$\begin{pmatrix} \mathrm{D}_{(p_2)^2}f(p_1) & \cdots & \mathrm{D}_{p_2 p_{n-1}}f(p_1) & \mathrm{D}_{p_2 p}f(p_1) \\ \vdots & \ddots & \vdots & \vdots \\ \mathrm{D}_{p_{n-1}p_2}f(p_1) & \cdots & \mathrm{D}_{(p_{n-1})^2}f(p_1) & \mathrm{D}_{p_{n-1}p}f(p_1) \\ \mathrm{D}_{pp_2}f(p_1) & \cdots & \mathrm{D}_{pp_{n-1}}f(p_1) & \mathrm{D}_{p^2}f(p_1) \end{pmatrix},$$

has multidegree $\bigl((n-1)\cdot(d-2), 2, \ldots, 2\bigr)$ in $p_1, \ldots, p_{n-1}, p$, so the argument of (3.21) produces a homogeneous polynomial $H$ of degree

$$(n-1)(d-2) + 2(d-2) = (n+1)(d-2)$$

in the variable $p_1$.

**(3.23)** The analysis of the divisibility of $\varphi$ by $T$ on the other hand will give rise to a couple-nodal homogeneous polynomial $K$ in the variable $p_1$.

Let us first compute the multidegree of $\varphi$ in the variables $p_1, \ldots, p_{n-1}, p$. This goes again as in (3.16). First of all, $\mathrm{Disc}(H_p \cap X)$ is $n$-homogeneous of $n$-degree $(d^\vee, \ldots, d^\vee)$, with $d^\vee = d(d-1)^{n-1}$, and so is its homogeneous piece of lowest degree with respect to the Zariski grading. Then $T$ has visibly $n$-degree $(d-1, 0, \ldots, 0, 1)$, while $\mathrm{Disc}(F_2)$ has $n$-degree $((n-1)(d-2), 2, \ldots, 2)$ as we saw in (3.22). The same computation gives the $n$-degree of $\mathrm{Disc}(\bar{F}_2)$, viz. $((n-2)(d-2), 2, \ldots, 2, 0)$. Eventually, one finds that $\varphi$ has degrees

$$d\bigl[(d-1)^{n-1} - 1\bigr] - 3(n-1)(d-2) \quad \text{and} \quad d(d-1)^{n-1} - 6$$

in $p_1$ and $p_2, \ldots, p_{n-1}, p$ respectively (note that the former degree is divisible by $d-2$). Therefore, we get by (3.21) a couple-nodal polynomial $K$ in $p_1$, homogeneous of degree

$$d\big[(d-1)^n - 1\big] - 3(n+1)(d-2) = (d-2)\Big(d \cdot \frac{(d-1)^n - 1}{d-2} - 3(n+1)\Big).$$

One thus obtains the following result.

**(3.24) Theorem.** *Let $X$ be a smooth degree $d$ hypersurface in $\mathbf{P}^n$, $n > 1$. The number of bitangent planes to $X$ passing through $n-2$ general points in $\mathbf{P}^n$ is*

$$\frac{1}{2}d(d-1)^{n-2}(d-2)\Big(d \cdot \frac{(d-1)^n - 1}{d-2} - 3(n+1)\Big).$$

For $n = 3$ one recovers Theorem (3.17), and for $n = 2$ the number of bitangents to a smooth plane curve of degree $d$, viz.

$$\frac{1}{2}(d^\vee - 1)(d^\vee - 2) - 3d(d-2) - \frac{1}{2}(d-1)(d-2) = \frac{1}{2}d(d-2)(d-3)(d+3).$$

# 4 – Two further enumerative computations

In this final section we present two enumerative computations for surfaces in $\mathbf{P}^3$, also taken from Salmon's book, which are close in spirit to the previous considerations.

## 4.1 – The flecnodal polynomial

This is carried out by Salmon in [16, §588], with [16, §473] as a fundamental tool. This has already been revisited in modern standards in [2], and actually extended there to hypersurfaces in a projective space of arbitrary dimension, so we are going to be brief.

**(4.1) The problem.** Let $S$ be a smooth surface in $\mathbf{P}^3$ of degree $d > 1$. For a general point $p \in S$, there are two lines having intersection multiplicity at $p$ with $S$ strictly greater than 2, namely the tangent lines to the two smooth branches at $p$ of the curve $\mathbf{T}_p S \cap S$, which intersect $S$ with multiplicity 3 at $p$. We shall see that those points $p \in S$ such that there is a line intersecting $S$ with multiplicity strictly greater than 3 at $p$ is a curve $F(S)$, cut out on $S$ by a polynomial of degree $11d - 24$. We call this curve (resp. polynomial) the *flecnodal* curve (resp. polynomial) of $S$.

At a general point $p$ of the flecnodal curve, the section of $S$ by its tangent hyperplane $\mathbf{T}_p S$ is a curve with a non-degenerate double point at $p$ (i.e., a double point with tangent cone of maximal rank), with one of its two local branches having an inflexion point at $p$. In general, the tangent line to the latter branch meets $S$ with multiplicity 4 at $p$. Those points $p \in S$ such that the curve $\mathbf{T}_p S \cap S$ has a tacnode (i.e., a double point with local equation $y^2 = x^4$) also belong to the flecnodal curve: they are its intersection points with the Hessian of $S$, and they are cuspidal points of the cuspidal double curve of $S^\vee$ (the latter curve parametrizes those hyperplanes that cut out a cuspidal curve on $S$).

The following statement is definitely a result in reduced elimination theory, although it does not strictly fit in the framework of Section 1.

**(4.2) Proposition** (see [16, §473]). *Let $f_q(p), g_q(p), h_q(p)$ be three bi-homogeneous polynomials in $p, q \in \mathbf{P}^3$ of bi-degrees $(\lambda, \mu)$, $(\lambda', \mu')$, $(\lambda'', \mu'')$ respectively. We assume that for the generic point $q \in \mathbf{P}^3$,*
$$\mathrm{mult}_q\big(V(f_q, g_q, h_q)\big) = \lambda\lambda'\lambda''.$$

*The locus of those $q \in \mathbf{P}^3$ such that $V(f_q, g_q, h_q)$ contains a point in addition to $q$ counted with multiplicity $\lambda\lambda'\lambda''$ is the zero locus of a homogeneous polynomial of degree*
$$\lambda'\lambda''\mu + \lambda\lambda''\mu' + \lambda\lambda'\mu'' - \lambda\lambda'\lambda''.$$

Of course the condition that the scheme $V(f_q, g_q, h_q)$ contains a point in addition to $q$ is equivalent to its having positive dimension. Salmon claims that it is equivalent to the fact that $V(f_q, g_q, h_q)$ contains a line; we were not able to prove this, but it turns out that this is not needed for the application.

*Proof.* We want to characterize when the scheme $V(f_q, g_q, h_q)$ has positive dimension. The idea is that this is equivalent to its having non-empty intersection with any hyperplane. So let $\ell$ be a non-zero linear form in $p$, and consider the resultant $\mathrm{Res}(\ell, f_q, g_q, h_q)$. It follows from the Poisson formula (see, e.g., [2, Prop. 2.2] and the references therein) and our assumption on $f, g, h$ that there exists a polynomial $R$ such that

(4.2.1) $$\mathrm{Res}(\ell, f_q, g_q, h_q) = \ell(q)^{\lambda\lambda'\lambda''} \cdot R(\ell, f_q, g_q, h_q).$$

Computing degrees, one sees that $R$ is homogeneous of degree 0 in the coefficients of $\ell$, i.e., it does not depend on $\ell$. It follows that $V(f_q, g_q, h_q)$ has positive dimension if and only if $R(f_q, g_q, h_q) = 0$. Eventually, one computes the degrees of $R(f_q, g_q, h_q)$ using the identity (4.2.1). □

**(4.3) Theorem.** *Let $S$ be a smooth surface in $\mathbf{P}^3$ of degree $d$. There exists a homogeneous polynomial $F$ of degree $11d - 24$ such that the locus of points $p \in S$ such that there is a line intersecting $S$ with multiplicity at least 4 in $p$ is cut out on $S$ by $V(F)$.*

*Proof.* Let $p \in S$. It follows from Theorem (3.5) that there is a line intersecting $S$ with multiplicity at least 4 in $p$ if and only if the three polar hypersurfaces $\mathrm{D}^1 S(p), \mathrm{D}^2 S(p), \mathrm{D}^3 S(p)$ have a common point besides $p$, and that this is equivalent to their having a whole line in common.

On the other hand, $\mathrm{D}^1 S(p), \mathrm{D}^2 S(p), \mathrm{D}^3 S(p)$, which have degrees $1, 2, 3$ respectively, intersect with multiplicity 6 at $p$ by Corollary (3.6). We are therefore in a position to apply Proposition (4.2), with $(\lambda, \mu), (\lambda', \mu'), (\lambda'', \mu'')$ equal to $(1, d-1), (2, d-2), (3, d-3)$ respectively. The result follows. □

### 4.2 – Number of bitangent lines

The Grassmannian of lines in $\mathbf{P}^3$ has dimension 4. Passing through a fixed point imposes 2 conditions to a line in $\mathbf{P}^3$, and being tangent to a surface (at an unprescribed point) imposes 1 condition, so one expects finitely many bitangent lines to a surface passing through a general point in $\mathbf{P}^3$. In this subsection we prove the following, along the lines of [15, §279].

**(4.4) Theorem.** *Let $S$ be a smooth surface of degree $d$ in $\mathbf{P}^3$, and $p \in \mathbf{P}^3$ a general point. The number of lines bitangent to $S$ and passing through $p$ is*

(4.4.1) $$\frac{1}{2}d(d-1)(d-2)(d-3).$$

Salmon's strategy is similar in spirit to that exposed in subsection 3.2. The fundamental fact is the following. Consider as before

$$\tilde{\mathcal{C}}_S = \{(p', p'') \in \mathbf{P}^3 \times \mathbf{P}^3 : p' \in S \text{ and } p'' \in \mathbf{T}_{p'}S - p'\}.$$

**(4.5) Proposition.** *There exists a polynomial $R$, bihomogeneous in $p', p''$ with respective degrees $(d-2)(d-3)$ and $(d+2)(d-3)$, such that the locus of points $(p', p'') \in \tilde{\mathcal{C}}_S$ such that the line $\langle p', p''\rangle$ is bitangent to $S$ is cut out on $\tilde{\mathcal{C}}_S$ by $V(R)$.*

*Proof.* Let $f$ be an equation of $S \subset \mathbf{P}^3$. For $(p', p'') \in \tilde{\mathcal{C}}_S$, we consider the homogeneous polynomial in $(\alpha, \beta)$

$$f(\alpha p' + \beta p'') = \alpha^d f(p') + \alpha^{d-1}\beta D_{p''} f(p') + \alpha^{d-2}\beta^2 D_{(p'')^2} f(p') + \cdots + \beta^d f(p'')$$
(4.5.1) $$= \beta^2 \left(\alpha^{d-2} D_{(p'')^2} f(p') + \cdots + \beta^{d-2} f(p'')\right)$$

(this is an abuse of notation: actually one should consider two liftings $\hat{p}', \hat{p}'' \in \mathbf{k}^4$ of $p'$ and $p''$ respectively, and $f(\alpha \hat{p}' + \beta \hat{p}'')$). Let $^\flat f(\alpha, \beta)$ be the homogeneous polynomial of degree $d-2$ between parentheses at the right-hand-side of (4.5.1). The line $\langle p', p''\rangle$ is a bitangent to $S$ if and only if $^\flat f$ has a multiple root, so the polynomial $R$ we are looking for is merely the discriminant of $^\flat f$, which is, up to sign, nothing but the reduced discriminant of $f$ as noticed in (2.11).

Now, the polynomial $f$ is of the form

$$f = \beta^2 \left(a_2 \alpha^{d-2} + a_3 \alpha^{d-3}\beta + \cdots a_d \beta^{d-2}\right)$$

where the coefficients $a_i$ are homogeneous in $p'$, respectively $p''$, of degree $d - i$, respectively $i$. Therefore, we deduce from (2.11) (see also (2.14)) that $\mathrm{redDisc}(f)$ has degrees $(d-2)(d-3)$ and $(d+2)(d-3)$ in $p'$ and $p''$ respectively, which concludes the proof. □

**(4.6) Remark.** The degree of $\mathrm{redDisc}(f)$ in $p''$ may be computed alternatively as follows. The plane curve $C_{p'} := \mathbf{T}_{p'}S \cap S$ has in general a double point at $p'$, and what we want is the number of lines in $\mathbf{T}_{p'}S$ passing through $p'$ and tangent to $C_{p'}$ at some other point. This is the number of ramification points of the projection of $C_{p'}$ from $p'$; the latter is a $(d-2):1$ map $\bar{C}_{p'} \to \mathbf{P}^1$, where $\bar{C}_{p'}$ denotes the normalization of $C_{p'}$ at $p'$, so it follows from the Riemann–Hurwitz formula that the number of ramification points equals

$$2g(\bar{C}_{p'}) - 2 + 2(d-2) = \left((d-1)(d-2) - 2\right) - 2 + 2(d-2) = (d+2)(d-3)$$

as required.

We were not able to find, on the other hand, a geometric argument to compute the degree in $p'$ in a comparable fashion. We wonder wether there is an explanation to why this degree $(d-2)(d-3)$ is so nice, in particular in its role in (4.4.1). It is conceivable that it has something to do with the degree of the dual to a smooth plane curve of degree $d-2$.

*Proof of Theorem (4.4).* Let $p \in \mathbf{P}^3$ be a general point. The locus of those points $p' \in S$ such that $(p', p) \in \tilde{\mathcal{C}}_S$ (equivalently, $p \in \mathbf{T}_{p'}S$) is the apparent boundary $D_p S \cap S$. Among these points $p'$, the locus of those $p'$ for which the line $\langle p, p'\rangle$ is a bitangent to $S$ is cut out by $V(R)$, where $R$ is the polynomial of Proposition (4.5); it is therefore a complete intersection in $\mathbf{P}^3$ of type $(d, d-1, (d-2)(d-3))$. One concludes by observing that there are two points $p'$ for every bitangent line to $S$ passing through $p$. □

Laurent Busé. Université Côte d'Azur, Inria, 2004 route des Lucioles, 06902 Sophia Antipolis, France. `laurent.buse@inria.fr`

Thomas Dedieu. Institut de Mathématiques de Toulouse ; UMR5219. Université de Toulouse ; CNRS. UPS IMT, F-31062 Toulouse Cedex 9, France. `thomas.dedieu@math.univ-toulouse.fr`

ThD was membre of project FOSICAV, which has received funding from the European Union's Horizon 2020 research and innovation programme under the Marie Skłodowska-Curie grant agreement No 652782.